\documentclass[11pt]{article}
\usepackage{graphicx}
\usepackage{amsmath,
amssymb,
fancybox,enumerate,color}

\newtheorem{theorem}{\bf Theorem}[section]
\newtheorem{lemma}{\bf Lemma}[section]
\newtheorem{proposition}{\bf Proposition}[section]
\newtheorem{corollary}{\bf Corollary}[section]
\newtheorem{remark}{\bf Remark}[section]

\newtheorem{definition}{\bf Definition}[section]

\newcommand{\ds}\displaystyle

\def\R{{\mathbb R}}

\makeatletter
 \@addtoreset{equation}{section}
\makeatother

\begin{document}
\title{Non-uniqueness for a critical heat equation \\  in two dimensions with singular data
}
\author{
Norisuke Ioku*\\
Graduate School of Science and Engineering, \\
Ehime University\\
Matsuyama, Ehime 790-8577, Japan,\vspace{5pt}\\
\quad\\
Bernhard Ruf
\\
Dipartimento di Matematica ``F. Enriques'',\\
Universit\`a degli Studi di Milano\\
via C. Saldini 50, Milano 20133, Italy,
\\
\quad\\
and\\
\quad\\
Elide Terraneo\\
Dipartimento di Matematica ``F. Enriques'',\\
Universit\`a degli Studi di Milano\\
via C. Saldini 50, Milano 20133, Italy.
}

\author{Norisuke Ioku, Bernhard Ruf and Elide Terraneo}
\date{
}
\maketitle



\abstract{Nonlinear heat equations in two dimensions with singular initial data are studied. In recent works nonlinearities with exponential growth of
Trudinger-Moser type have been shown to manifest critical behavior: well-posedness in the subcritical case and
non-existence for certain supercritical data.
In this article we propose a specific model nonlinearity with
Trudinger-Moser growth  for which we obtain surprisingly complete
results: a) for initial data strictly below a certain singular
threshold function $\widetilde u$ the problem is well-posed, b)
for initial data above this threshold function $\widetilde u$,
there exists no solution, c) for the singular initial datum
$\widetilde u$ there is non-uniqueness.  The function $\widetilde u$ is a weak
stationary singular solution of the problem, and we show that
there exists also a regularizing classical solution with the same
initial datum $\widetilde u$. }

\section{Introduction}
Consider  the following Cauchy problem with Dirichlet boundary condition
\begin{equation}\label{eq1}
 \left\{
  \begin{split}
   &\partial_tu-\Delta u
    =f(u)
    & \text{in}\ \ \Omega, \ t>0,
    \\
   &u(t,x)
    =0
  \ \   & \ \ \text{on}\ \partial \Omega, \ t>0,\\
    &u(0,x)=u_0(x) &\text{in}\ \ \Omega,
  \end{split}
 \right.
\end{equation}
where $\Omega$ is an open domain in $\mathbb R^N$. It is
well-known that for bounded initial data $u_0$ and for
$C^1$-nonlinearities $f$, this equation has a local-in-time
solution $u \in L^{\infty}_{loc}((0,T];L^{\infty}(\Omega))$ for
some $T>0$. In this article we address some questions concerning
singular initial data $u_0 \notin L^\infty(\Omega)$.
The case of power-type nonlinearity $f(s)=|s|^{p-1}s$ has been
widely studied beginning with the seminal works of F. Weissler
(see \cite{BC, NS, Te, W0, W1, W2, W3} and Section \ref{2} for a
description of known results). Let us focus our attention to the
so-called critical   nonlinearity $f(s)=|s|^{\frac 2{N-2}}s$,
($N\geq 3$) and let us  consider initial data in the Lebesgue
space $L^{\frac{N}{N-2}}(\mathbb{R}^N)$, which is invariant under
the scaling of the equation   and which  has the same
integrability as the growth of the nonlinearity. In this case the
existence and uniqueness of a local-in-time (classical) solution for
any initial data  hold.  However, some non-uniqueness phenomena of
(distributional) solutions appear. Moreover, for small data the
solution exists globally in time.

 In dimension $N=2$ this case does not happen and  one may expect a critical  situation for certain nonlinearities with higher than
polynomial growth. In recent works \cite{IJMS,I,IRT,RT,FKRT}
(see also \cite{MT} for more general
nonlinearities)
it was shown that nonlinearities with Trudinger-Moser growth, see
\cite{Po,Tr,M},
\begin{equation}\label{exp}
f(s) \sim e^{s^2} \  \ \hbox{ for } \ |s| \ \hbox{ large}
\end{equation}
in conjunction with data from the Orlicz space
$$
\exp L^2(\Omega) := \big\{ u
\in L_{loc}^1(\Omega) : \int_\Omega\left( e^{\alpha u^2}-1\right)
dx < \infty \ \  \hbox{for some } \alpha > 0\big\},
$$
show some of the critical behavior (see Section \ref{2}, Remark \ref{dc}):
\par \smallskip \noindent
-- existence of global-in-time solutions for small data $u_0$ in $\exp L^2(\Omega)$;
\par \smallskip \noindent
-- non-existence of solutions for some large initial data $u_0 \in \exp L^2(\Omega)$;
\par \noindent
-- existence of local-in-time solutions for any initial data $u_0 \in \exp L_0^2(\Omega) :=
\overline{C_0^\infty(\Omega)}^{\|\cdot\|_{\exp L^2}}$.
\par \medskip
In this paper we set out to complete the picture by proving a
non-uniqueness result for a particular equation on a ball
$B_\rho(0) \subset \mathbb R^2$. Indeed, for a certain
nonlinearity
$f(t)$
with growth of type \eqref{exp}
(more precisely, see \eqref{eq:ioku1.2})
 we show the existence of a singular solution $\widetilde u \in \exp L^2(B_\rho)$ for the corresponding elliptic equation, which  gives rise to a singular stationary distributional solution of the parabolic equation. The solution $\widetilde u$ has the asymptotic profile $\widetilde u(x) \sim \sqrt{-2\log|x|}$, for $|x|$ small, and belongs to $\exp L^2 \setminus \exp L^2_0$. We prove furthermore that the same initial datum $\widetilde u$ gives also rise to a regularizing solution, and hence we have  non-uniqueness.
\par \smallskip \noindent
Indeed, for this particular initial datum $\widetilde u$ and the nonlinearity
$f(t)$,
we get the following surprisingly complete result:
\par \medskip \noindent
 {\bf Theorem A} \label{theorem:1.1}
\ {\it Let the initial datum $u_0$ for the problem
\begin{equation}\label{g} \partial_t u - \Delta u =
f(u)
  \ \hbox{ in } \ B_\rho(0) \ , \
  u = 0 \ \hbox{ on } \  \partial B_{\rho}(0),
\end{equation} be given by $u_0(x)=\mu\, \widetilde u(x)$, $ \mu > 0$. Then the following hold:
\par \smallskip \noindent
1) {\rm (well-posedness)} \ If \ $\mu < 1$, then the equation has a unique regular local-in-time solution.
\par \smallskip \noindent
2) {\rm (non-uniqueness)}\ If  \ $\mu = 1$, then $u_0 = \widetilde u$ is a singular (distributional) stationary solution, and there exists a regular solution with the same initial datum $\widetilde u$.
\par \smallskip \noindent
3) {\rm (non-existence)} \ If  \ $\mu > 1$, then the equation has no non-negative solution, in any positive time
interval.
}

\par \vspace{0.6cm}

In Section 2 we present more detailed motivations and some background for this problem, and a more precise statement of our results. We point out that the phenomena described in Theorem A are rather subtle, and the function spaces (Orlicz and Lorentz spaces) and related notions of solution have to be chosen very carefully. After introducing these concepts, we formulate a precise statement of Theorem A in Theorem \ref{th2.1}, see end of Section 2.
\par \smallskip \noindent
In Section 3 we give some preliminary results on the heat kernel in Orlicz spaces and Lorentz spaces which will be needed in the proofs, and the notions of solution (weak, classical) will be introduced.
\par \smallskip \noindent
In Section 4 we construct a singular solution $\widetilde u(x)$ of the elliptic equation \eqref{g}: we use that
$\sqrt{-2\log|x|}$ is an exact solution of \eqref{g} for large values of  $\widetilde u(x)$, and then employ the shooting method to construct a solution with zero boundary values on a suitable ball $B_\rho$.
\par \smallskip \noindent
In Section~5 we prove the well-posedness of equation \eqref{g} for initial data below the threshold function $\widetilde u$, i.e. statement 1 in Theorem A and Theorem \ref{th2.1}. This is done with a contraction argument in a suitable function space.
\par \smallskip \noindent
In Section~6 we prove the non-uniqueness result (statement 2 of Theorem A  and of Theorem \ref{th2.1} below).
The stationary singular solution is given by $\widetilde u(x)$, as obtained in Section 4. The existence of a regular solution
with the same initial datum $\widetilde u(x)$ is quite delicate: we first consider an auxiliary equation in a Lorentz space
setting with a cubic nonlinearity and with initial datum which belongs to the Lorentz spaces $L^{2,q}$ for all $q > 2$, but not for $q = 2$.
 From this solution we then produce, by a suitable transformation (inspired by Brezis-Cazenave-Martel-Ramiandrisoa
 \cite{BCMR} and Fujishima-Ioku \cite{FI}), a super-solution of the Cauchy problem \eqref{g}. Finally, applying Perron's monotone method,
 we then obtain a classical solution of problem \eqref{g}.
\par \smallskip \noindent
In Section~7 we give the proof of the non-existence result (statement 3 in Theorem A and Theorem \ref{th2.1}).
We show that for data above the threshold function $\widetilde u(x)$ we encounter instantaneous blow-up, i.e. for no
positive time $T$ can a solution exist.
\par \medskip
We expect that similar phenomena hold in more general situations, but we note that the growth of the nonlinearity,
 the behavior of the singular initial data, and the employed function spaces will have to be very carefully calibrated.

\section{Origin of the problem and main result}\label{2}
\subsection{Polynomial nonlinearities}
The study of equation \eqref{eq1} with singular
data began with the pioneering works of F. Weissler \cite{W0},
\cite{W1}.
 He
considered equation \eqref{eq1} on the whole space $\mathbb
R^N$, with power type nonlinearities $f(s) = |s|^{p-1}s$ and with
singular data in certain Lebesgue spaces $L^q(\mathbb R^N)$. For
power nonlinearities the equation \eqref{eq1} enjoys a scale
invariance: if $u$ is a solution, then also
$$u_\lambda(t,x) := \lambda^{2/(p-1)}\, u(\lambda^2 t,\lambda x)
$$
is a solution.
One notes that the initial data space $L^q(\mathbb R^N)$ is invariant under this scaling if and only if $q = q_c = \frac{N(p-1)}2$.
  This exponent serves as a limiting or {\it critical} exponent for the well-posedness of
  the Cauchy problem \eqref{eq1} with  $f(s) = |s|^{p-1}s$ and initial data $u_0\in L^q(\mathbb
  R^N)$.
Indeed one has:
\begin{itemize}
\item[--]   if $q > q_c$, $q\geq 1$ or $q = q_c$,
$q>1$, then  the Cauchy problem \eqref{eq1} has a unique
local-in-time solution in $C([0,T], L^q(\mathbb{R}^N))\cap
L^{\infty}_{loc}(
(0,T),
L^{\infty}(\mathbb{R}^N))$ for some $T>0$,
(see \cite{BC}, \cite{W0},  \cite{W1}).
Moreover, in the critical case $q=q_c,  q>1$, for
sufficiently small data in $L^{q_c}(\mathbb R^N)$ there exist
 {\it global-in-time} solutions  (see \cite{W2});


\item[--]  if $1 \le q < q_c$, then there exist  some non-negative
initial data in $L^q(\mathbb R^N)$ for which there is no
non-negative solution for any positive time $T > 0$ (see
\cite{BC},\cite{W1}, \cite{W3}).

\end{itemize}

For $q\geq p$,  then $C([0,T],
L^q(\mathbb{R}^N))\subset L^{p}_{loc}((0,T)\times\mathbb{R}^N))$
and for any $u\in C([0,T], L^q(\mathbb{R}^N))$ each term of
equation \eqref{eq1} is a distribution. Therefore for $q\geq q_c$,
$p\geq q$, one may ask whether the  solution obtained by Weissler
is unique in the larger class $C([0,T], L^{q}(\mathbb{R}^N))$. The
known results are:
\begin{itemize}
\item[--]
 if $q>q_c$, $q\geq  p$ or  $q=q_c$,  $q>p$ uniqueness still holds in the class  $C([0,T], L^{q}(\mathbb{R}^N))$ (see \cite{BC}, \cite{W0}).
\end{itemize}
In the  case  $q = q_c$ and $q = p$, then $q=p = \frac N{N-2}$ which is referred to as  {\it
doubly critical case}  in \cite[Remark 5]{BC},
 Ni-Sacks \cite{NS} proved  that (for the unit ball $B_1 \subset \mathbb R^N$) there  exists a
 {\it stationary  singular} solution -- which  is different from the regularizing solution of Weissler.
  This non-uniqueness result was extended to the whole space $\mathbb R^N$ by Terraneo 
\cite{Te}.

 We remark that  if $p>\frac N{N-2}$  there exists an
explicit singular stationary solution of \eqref{eq1} with $f(s) =
|s|^{p-1}s$ in $\mathbb{R}^N$. This is another way in which
$p=\frac N{N-2}$ is critical and so we can say that  $q = q_c =
\frac N2(p-1) = p$, i.e. $p = \frac N{N-2}$, is doubly critical.

  \begin{remark}\label{dc} Note that the "doubly critical"
case is characterized by the simultaneous appearence of the following two phenomena:
\begin{itemize}
\item[-] global-in-time existence for small data; \item[-]
non-uniqueness for some data.
\end{itemize}
\end{remark}
\par \smallskip

\subsection{The limiting case: the $ H^s$ - $ L^p$ correspondence}
Note that in $\mathbb R^2$ the double critical exponent $q=q_c =p=
\frac N{N-2}$ becomes infinite. If we look for a suitable
"critical growth" in two dimensions, we may be
  guided by recent results for dispersive equations.
\par \smallskip
Indeed, for the corresponding Nonlinear Schr\"odinger equation,
where one works with energy methods,
 one has similar phenomena for initial data $u_0$ in Sobolev spaces $H^s(\mathbb R^N)$: again one finds,
 corresponding to the power nonlinearity $|u|^{p-1}u$, an associated critical space $H^{s_c}(\mathbb R^N)$ with
 $s_c = \frac N2 - \frac 2{p-1}$.
 Cazenave-Weissler \cite{CW} showed local-in-time existence for all $u_0 \in H^s(\mathbb R^N)$ for $s \ge s_c$, and global-in-time existence for
  small data for $s = s_c$.  The critical exponents for the $H^s$-theory for the Nonlinear Schr\"odinger and  heat equations  coincide,
   while the critical exponents for the $H^s$-theory and the  $L^p$-theory for the heat equation are related by the Sobolev
   embedding:
   $H^{s_c} \subset L^{q_c}$, with $q_c = \frac
{2N}{N-2s_c} = \frac N2 (p-1)$.

In the  limiting critical case $s_c = \frac N2$ we have again that
$H^{N/2}(\mathbb R^N) \subset L^q(\mathbb R^N)$, for all $q \ge
1$, but $H^{N/2}(\mathbb R^N) \not \subset L^\infty(\mathbb R^N)$.
By a result by S. Pohozaev \cite{Po} and N. Trudinger \cite{Tr} we
know that for $u \in H^{N/2}$ one has $\int_{\mathbb R^N}
(e^{u^2}-1) dx < \infty$, and this is the maximal growth for
integrability. Using nonlinearities with this type of growth in
the Nonlinear Schr\"odinger equation (NLS equation)
$$
i\partial_t u + \Delta u = f(u) \ \ \hbox{ with } \ f(u) \sim e^{u^2}
$$
Nakamura-Ozawa \cite{NO} were indeed able to prove a global-in-time existence result for small initial data in $H^{N/2}(\mathbb R^N)$,  and so in particular in $H^1(\mathbb R^2)$  for
$N = 2$.
For other related results we refer to \cite{CIMM}.
\par \medskip
\subsection{Back to the heat equation}
The result of Nakamura-Ozawa was recently transposed to the heat equation by Ibrahim-Jrad-Majdoub-Saanouni \cite{IJMS}, showing  {\it local}-in-time existence and uniqueness for the equation \eqref{eq1}, with $f(u) \sim e^{u^2}$, $x \in \mathbb R^2$, and for any initial data $u_0 \in H^{1}(\mathbb R^2)$. Two observations are in order:
\vspace{-0.1cm}
\begin{itemize}
\item[--]  the initial data space $H^1(\mathbb R^2)$ is natural for the NLS equation,
 where one works with energy methods, but less so for the heat equation,
 where an integrability condition on the initial data ought to be
 sufficient;
\vspace{-0.1cm}
\item[--] by Nakamura-Ozawa \cite{NO} a {\it global-in-time} result holds
for the NLS equation with $f(u)\sim e^{u^2}$, for {\it small data} in $H^1(\mathbb R^2)$; comparing with the critical case for polynomial nonlinearities, one can say that $f(u)\sim e^{u^2}$ behaves like a critical growth nonlinearity for the  NLS equation. However, the uniqueness result in \cite{IJMS} suggests that $f(u)\sim e^{u^2}$ with initial data in $H^1(\mathbb{R}^2)$ is not a double critical case (in the sense of  Remark \ref{dc}).
\end{itemize}

  Here we are looking, in dimension $N=2$, for a
data space which has similar "double critical" phenomena as
described in Remark \ref{dc}.  We propose the Orlicz space
determined by the mentioned estimates by Pohozaev and Trudinger,
namely $H^1(\mathbb{R}^2) \subset L^{\varphi}(\mathbb{R}^2)$ with
Young-function $\varphi(t) = e^{t^2} -1$ (for details, see Section \ref{O}  below). We will denote this space by
  $\exp L^2(\mathbb R^2) := L^\varphi(\mathbb R^2)$.
 In fact, in \cite{RT,I,IRT},
 small-data global-existence and large-data non-existence result were proved for this space.
%
%
\par \smallskip
In this paper, we focus on the following particular case of an exponential nonlinearity with Trudinger-Moser growth.
Consider the nonlinearity $f(s)$ given by
\begin{equation}\label{eq:ioku1.2}
f(s) := \left\{\begin{array}{ll}\ds  \frac 1 {|s|^3} \ e^{s^2} \  \ &\hbox{ if } \ |s| > \beta, \vspace{0.2cm} \\
\ds \alpha s^2 \ &\hbox{ if } \  |s| \le \beta
\end{array}\right.
\end{equation}
with  $\alpha=\frac{{\rm e}^{5/2}}{\left(5/2\right)^{5/2}}$  and
$\beta =\sqrt{\frac 52}$  such that the function $f $ belongs to
$C^1(\mathbb{R})$, it is  increasing on $[0,+\infty)$ and  convex
on $\mathbb{R}$. We will show that the nonlinearity
\eqref{eq:ioku1.2}, together with suitable initial data,
shows all the phenomena of a double critical case  for the 2-dimensional problem,
with respect to  existence, non-existence, uniqueness and
non-uniqueness.
\par \smallskip
To this end, we first prove the existence of a radial singular
solution for the Dirichlet boundary value problem in $B_{\rho}
\subset \mathbb R^2$
\begin{equation}\label{eq:1.2}
 \left\{
  \begin{split}
   -\Delta u
    &=f(u)
    & \text{in}\ \ B_{\rho},
    \\
   u(x)
    &=0
    & \text{on}\ \partial B_{\rho}
  \end{split}
 \right.
\end{equation}
for some $\rho>0$. By a {\it singular} solution we mean a
 solution which belongs to
$C^2(B_{\rho}~\setminus~\{0\})$,
  which is unbounded on $B_{\rho}$ and which satisfies the elliptic equation in the sense of distributions on $B_{\rho}$.
 Moreover  this solution $\widetilde u$ belongs to the Orlicz space $\exp L^2(B_\rho)$. More precisely, we prove
\par \medskip
\begin{proposition}\label{prop1}
There exist a constant $\rho>0$ and a function $\widetilde u\in
C^2(B_{\rho}\setminus\{0\})\cap
C(\overline{B}_{\rho}~\setminus~\{0\})$ which is a classical
solution on $B_{\rho}\setminus\{0\}$  for the Dirichlet boundary
value problem  \eqref{eq:1.2}. Moreover, the following hold:
\begin{itemize}
\item[\rm (i)]
$\widetilde u(x)= \sqrt{-2\log(|x|)} $ \ in a neighborhood
of the origin;
\item[\rm (ii)] $\widetilde u$ is a solution of the elliptic equation
\eqref{eq:1.2} on $B_{\rho}$ in the sense of distributions.
\end{itemize}
\end{proposition}

 \begin{remark} $ $ \par \smallskip \noindent
a) With the change of variable $y = \frac x \rho$  and the corresponding changes in the nonlinearity
$ f(u) \rightsquigarrow \rho^2 f(u)$ and initial datum\ $\widetilde u(x) \rightsquigarrow \widetilde u(\frac x \rho)$\ the equation can be considered on $B_1(0) \subset \mathbb R^2$.
 \par \smallskip \noindent
 b) The nonlinearity $f(s)$ may be generalized to
 $$
 f(s) = \left\{ \begin{array}{lll} &\frac 1{|s|^3} \, e^{s^2} \quad &, \quad |s| > \beta_p \vspace{0.2cm} \\
 &\alpha_p\, s^p &, \quad |s| \le \beta_p
 \end{array} \right.
 $$
 for any choice of $p > 1$ and suitable values $\alpha_p, \beta_p$ (which are uniquely dependent on $p$ since $f(s)$ is required to be of class $C^1(\mathbb R)$).
\end{remark}

The particular form of the nonlinearity \eqref{eq:ioku1.2} is due to the existence of the (almost explicit) singular solution given in Proposition \ref{prop1}.(i). It would be of interest to prove the existence of singular distributional solutions for equation \eqref{eq:1.2} for more general nonlinearities.

\par \medskip
\subsection{Main result: A heat equation in 2-dimensions with  double critical phenomena}
Let us now  consider the
 following Cauchy problem with Dirichlet boundary condition on $B_{\rho}\subset \mathbb R^2$
\begin{equation}\label{eq:2.1}
 \left\{
  \begin{split}
   &\partial_tu-\Delta u
    =f(u)
    & \text{in}\ \ B_{\rho}, \ t>0,
    \\
   &u(t,x)
    =0
    & \quad \text{on}\ \partial B_{\rho}, \ t>0,\\
    &u(0,x)=u_0(x) &\text{in}\ \ B_{\rho},
  \end{split}
 \right.
\end{equation}
\par \smallskip \noindent
where  the nonlinear term $f(u)$ is defined in \eqref{eq:ioku1.2}. We will show that the singular function $\widetilde u$
obtained in Proposition \ref{prop1} yields a neat  separation into the cases of well-posedness,
 non-uniqueness and non-existence, and so we may say that we are in a ``double critical'' situation
 in the sense of Remark~\ref{dc}.

To state the theorem, we denote the Schwarz symmetrization of a
measurable function $\varphi : B_\rho\to \mathbb R$   by
$\varphi^\sharp$ (for details, see Section \ref{L}).
Moreover we introduce the complete metric space for $T,\
{\mu}^*>0$,
\begin{equation}
\label{eq:2.4aaa} M_{T,\,\mu^*}=\Big\{ u\in
L^{\infty}(0,T;\exp L^2(B_{\rho})): \ \ \sup_{t\in
(0,T)}\|u(t)\|_{L^f_\gamma(B_{\rho})}\leq \mu^*\Big\},
\end{equation}
 where $\|\cdot\|_{L^f_{\gamma}}$
is the Luxemburg norm defined by
\[
\|u\|_{L^f_\gamma(B)}=\inf
\left\{
\lambda>0: \ \int_Bf
\left(\frac
{|u(x)|}{\lambda}
\right) dx \leq \gamma
\right\}
\]
with
$\gamma = \int_{B_{\rho}}f(\tilde u(x))dx<\infty$.
For the
definitions of the Orlicz
space
$\exp L^2(B_{\rho})$
with the Luxemburg norm
$\|\cdot\|_{L_\gamma^f}$ under specific choice of $\gamma$,
and of weak
and $\exp L^2-$classical solutions, see Sections
 \ref{O} and \ref{3.4}.
%

%
%

\begin{theorem}\label{th2.1}
Let $\widetilde u$ denote the singular solution of the elliptic equation \eqref{eq:1.2} given by Proposition \ref{prop1}.
\par \medskip \noindent
1) (well-posedness) If the initial datum $u_0$ in \eqref{eq:2.1}
satisfies
\begin{equation}\label{mu}
\mu := \sup_{x\in B_{\rho}}\frac{{u_0}^{\sharp}(x)}{\widetilde u(x)}< 1,
\end{equation}
then problem \eqref{eq:2.1} is well-posed, i.e.  for any
 $\mu< \mu_1<1$ there exist a positive time
 $T=T(\mu_1)>0$ and
 a unique function\ $u$\ in
the complete metric space $M_{T,\, \mu_1}$
which is a weak solution of the Cauchy problem \eqref{eq:2.1}.
\noindent Furthermore,
it is an  $\exp L^2-$classical solution of \eqref{eq:2.1} on
$(0,T)\times B_\rho$.

\par \medskip \noindent
2) (existence and non-uniqueness) If the initial datum $u_0$ satisfies
\begin{equation}\label{mu1}
\mu = \sup_{x\in B_{\rho}}\frac{{u_0}^{\sharp}(x)}{\widetilde u(x)} \le 1,
\end{equation}
then  \eqref{eq:2.1}  admits an $\exp L^2-$classical solution $u$ in
some time interval $(0,T)$.  If $\mu<1$ this solution
belongs to $M_{T,\, \mu_1}$ for some $\mu< \mu_1 < 1$ (for
sufficiently small $T$), and hence coincides with the solution
obtained in 1). If $\mu = 1$  for any  $ \mu_2 > 1$ the solution belongs to $M_{T,\,
\mu_2}$ for some $ T$ and may not be unique in this space.

Indeed, for
$u_0 = \widetilde u$
the equation \eqref{eq:2.1} has, in addition to this classical
solution, the singular stationary (distributional) solution
$\widetilde u$  which belongs to $M_{T,1} \subset
M_{T,\,\mu_2}$.
\par \medskip \noindent

\par \medskip \noindent
3) (non-existence) Let $u_0 = \mu\,  \widetilde u$, with $\mu >
1$. Then the problem \eqref{eq:2.1} does not possess non-negative
$\exp L^2-$classical solutions on any positive time interval
$(0,T)$.
\end{theorem}

\begin{remark}\label{Remark:2.1aaa}$ $
\par \smallskip \noindent
 a) The solution in Theorem 2.1.1) can be continued
as long as $\mu(u(t)) := \sup_{x \in B_\rho} \frac {u^\sharp
(t,x)}{\widetilde u(x)} < 1$. If $\mu(u(t^*)) = 1$ for some $t^* >
0$, then the local theory fails and non-uniqueness may occur.
\par \smallskip \noindent
b) Since $\widetilde u$ is a radially symmetric and non-increasing function, the Schwarz symmetrization of $\widetilde{u}$
coincides with $\widetilde{u}$.
Therefore, Theorem~A 1) and 2) are particular cases of Theorem~\ref{th2.1} with
$u_0=\mu \widetilde{u}$, $0<\mu <1$ and $u_0=\widetilde{u}$, respectively.
\end{remark}
\begin{remark}
\par \smallskip \noindent
 We mention  that, with different techniques,  Galaktionov-Vazquez \cite{GV} and
Souplet-Weissler \cite{SW} proved similar results for the heat
equation with polynomial nonlinearity. Indeed, if $N>2$ and
$p>\frac N{N-2}$ the function
$V(x)=\beta^{1/(p-1)}|x|^{-2/(p-1)}$, where $\beta=\frac
2{p-1}\left(N-2-\frac2{p-1}\right)$ is an explicit stationary
distributional solution for the equation \eqref{eq1} with
$f(s)=|s|^{p-1}s$. For $N>2$ and for  any $\frac N{N-2}<p<p^*$
(where $p^*=+\infty$ if $N\leq 10 $ and
$p^*=\frac{N-2\sqrt{N-1}}{N-4-2\sqrt{N-1}}$  if $N>10$) the
equation  with initial data $\mu V(x)$, with $\mu\in
[1,1+\varepsilon)$ for $\varepsilon >0$ small enough, admits at
least a nonnegative regular solution $u(t)$ that converges to $\mu
V(x)$ in the sense of distributions as $t\to 0$. This implies similar phenomena
of non-uniqueness as in part 2) of Theorem \ref{th2.1}. Moreover,
for large values of $\mu$ the Cauchy problem with initial data
$\mu V(x)$ has no local nonnegative solution (see \cite{W3}).
\end{remark}

\section{Preliminary results}\label{sec3}
Let $B\subset \mathbb{R}^2$ be a ball centered at the origin. In
this section we recall some properties of Orlicz and Lorentz
spaces on $B$, and of the heat kernel in these spaces. We also
introduce  the definition of weak and $\exp L^2-$classical solution
of the problem \eqref{eq:2.1}.
\subsection{Orlicz spaces}\label{O}
Let us recall the definition of the Orlicz space $L^\varphi(B)$,
where $\varphi(u)$ is a Young function (convex, $\varphi(0) = 0$).
First we introduce the {\it Orlicz class} $K^\varphi(B)$ by
$$
K^\varphi(B) = \Big\{ u\in L^1(B): \int_B\varphi\Big(
{|u(x)|}\Big) dx <+\infty\Big\}.
$$
Then the {\it Orlicz space} $L^\varphi(B)$ is given by the linear
hull of the Orlicz class $K^\varphi(B)$
and its norm is given by the Luxemburg type
\[
  \|u\|_{L^{\varphi}(B)}
  :=
    \inf
    \left\{
     \lambda>0:
     \int_{B}
      \varphi\left(
       \frac{|u(x)|}{\lambda}
      \right)dx
       \le 1
    \right\}.
\]
For $\varphi(u) = e^{u^2}
-1$ we define $\exp L^2(B)=L^\varphi(B)$.
Let now $f$ be the convex function defined in \eqref{eq:ioku1.2}.
 Since for any $0< b <1$ there exist $C_1,C_2>0$
such that
 \begin{equation}\label{relor}
 C_1\left({\rm e}^{b\, u^2}-1\right)\leq f(u)\leq C_2\left({\rm
 e}^{u^2}-1\right),
 \end{equation} we have that the Orlicz space $\exp L^2(B)$ coincides with the Orlicz
  space generated by the convex function $f$, namely,
$$
\exp L^2(B)= L^f(B)
$$
and this space can be endowed  with the following equivalent
norm
\begin{equation}\label{Lf1}
\|u\|_{L^f_\gamma(B)}=\inf
\left\{
\lambda>0: \ \int_Bf
\left(\frac
{|u(x)|}{\lambda}
\right) dx \leq \gamma
\right\}
\end{equation}
for any fixed positive constant $\gamma$. Indeed, we have

\begin{proposition} \label{proposition:2.1} Let $\gamma>0$. There exist two positive
constants $c$, $C$ such that
\begin{equation}\label{Orl1}
c\| u\|_{L^f(B)}\leq \| u\|_{L^f_\gamma(B)}\leq C\| u\|_{L^f(B)}
\end{equation}
and
\begin{equation}\label{Orl2}
c\| u\|_{\exp L^2(B)}\leq \| u\|_{L^f(B)}\leq C\| u\|_{\exp L^2(B)}.
\end{equation}
Furthermore, in \eqref{Orl1} one may choose $c=\min(1,\frac
1\gamma)$ and $C=\max(1,\frac 1\gamma)$.
\end{proposition}

\noindent {\bf Proof.} Let us prove the first inequality. Assume
$0<\gamma<1$. By the definition
 we get directly $\|u\|_{L^f(B)}\leq \|u\|_{L^f_\gamma(B)}$.
 On the other hand
 thanks to the convexity of $f$ and the property $f(0)=0$ we
 obtain
 $$
 f\Big(\gamma \frac {u}\lambda\Big)= f\Big(\gamma \frac {u}\lambda+(1-\gamma)0\Big)\leq \gamma\, f\Big(\frac
 {u}\lambda\Big)+(1-\gamma) f(0)=\gamma \, f\Big(\frac
 {u}\lambda\Big).
 $$
Therefore  it holds
 \begin{equation*}
 \begin{aligned}
\|u\|_{L^f(B)}&=\inf \Big\{\lambda >0: \int_B\gamma f\Big(\frac
{u}\lambda\Big) dx\leq \gamma\Big\}\\
& \geq\inf \Big\{\lambda >0: \int_B f\Big(\frac {\gamma
{u}}\lambda\Big) dx\leq \gamma \Big\}\\
& =\gamma \|u\|_{L^f_\gamma(B)}.
\end{aligned}
 \end{equation*}
For $\gamma >1$ we can apply similar arguments to $0<\frac
1\gamma<1$. The second inequality follows from the relation
(\ref{relor}) and from the definition of Orlicz space (see
\cite[Section~8.4 and 8.12]{AF}).
This completes the proof of Proposition~\ref{proposition:2.1}.

In this paper we choose $\gamma:=\int_{B_{\rho}}f(\widetilde u(x))dx$.
It will be proved in Section~\ref{sec:sing} that $f(\widetilde u)$ is integrable, therefore $\gamma$ is well-defined.
This special choice of $\gamma$ is one of the keys to reach  a neat classification as in Theorem~\ref{th2.1}.

\subsection{Heat kernel}

Now we collect some  results concerning the solution of the heat
equation on the ball (see Appendix B in \cite{QS}). Let us denote
by $e^{t\Delta}$ the Dirichlet heat semigroup in $B$. It is known
that for any $\phi \in L^p(B)$, $1\leq p\leq +\infty$, the function
$u=e^{t\Delta}\phi$ solves the heat equation $u_t-\Delta u=0$ in $
(0,+\infty)\times B$ and $u\in C((0,+\infty)\times \overline{B})$,
$u=0$ on $ (0,+\infty)\times \partial B$. Moreover, there exists a
positive
${C}^\infty$
 function $G_B: B\times
B\times(0,+\infty) \to \R$ (the Dirichlet heat kernel) such that
$$
{\rm e}^{t\Delta }\phi(x)=\int _{B}G_B(x,y,t) \phi(y)dy,
$$
for any $\phi \in L^p(B)$, $1\leq p\leq +\infty$. We prepare
several basic lemmas.

\begin{lemma}\label{jensen}
Let $\phi: B \to [0,\infty)$ be a measurable function and
$H:\mathbb{R} \to \mathbb{R}$ be a convex function  such that
$H(0)=0$. Then
$$
H\left({\rm e}^{t\Delta }\phi\right)\leq{\rm e}^{t\Delta
}H\left(\phi\right).
$$
\end{lemma}

\noindent {\bf Proof.} Let $H$ be a convex function and $\phi\geq
0$ be a measurable function.  By  Jensen's inequality, denoting
$\overline{ G}=\overline{ G}(x,t)=\int _{B}G_{B}(x,y,t) dy$,
 we
obtain
\begin{equation*}
H\Big(\frac1 { \overline G(x,t)}\int_{B}G_B(x,y,t) {\phi(y)}
{dy}\Big)\leq \frac1 { \overline G(x,t)}\int_{B}G_B(x,y,t)H\big(
{\phi(y)}\big){dy}.
\end{equation*}
Therefore
\begin{equation}\label{I}
 H\Big(\frac{{\rm e}^{t\Delta }\phi}{\overline G
}\Big)\leq\frac1 { \overline G}\ {\rm e}^{t\Delta
}H\left({\phi}\right).
\end{equation} Moreover  by the  convexity of $H$,   the property
$H(0)=0$,  and $\overline G(x,t)\leq 1$ for any $x\in B$ and
$t>0$ we have
$$
H(s)= H\Big(\overline G \frac{s}{\overline G}+(1- \bar G)0 \Big)
\leq\overline G\ H\left(\frac{s}{\overline G}\right)
$$
and so  for $s={\rm e}^{t\Delta }\phi$ we get
\begin{equation}\label{II}
\frac{H(e^{t\Delta}\phi)}{\overline G}\leq H\Big(\frac{{\rm
e}^{t\Delta }\phi}{\overline G  }\Big).
\end{equation}
Finally, \eqref{I} and \eqref{II} imply  the desired inequality
$$
H\left({\rm e}^{t\Delta }\phi\right)\leq {\rm e}^{t\Delta
}H\left({\phi}\right).
$$
\begin{lemma}\label{heatcont}
There holds
$$
\|e^{t\Delta}\phi\|_{L^f_\gamma}\leq  \|\phi\|_{L^f_\gamma}
$$
for all $t>0$ and $\phi \in L^f_\gamma$(B).
\end{lemma}

\noindent {\bf Proof.} Here $f$ is the function in
\eqref{eq:ioku1.2}. Since $f$ is convex on $\mathbb{R}$ and
$f(0)=0$, it follows from the previous Lemma   and the property
$\overline{ G}(x,t)\leq 1$ for any $x\in B$ and $t>0$ that
\[
 \int_{B}f\Big(\frac{|e^{t\Delta}\phi|}{\lambda}\Big)dx
\le
 \int_{B}f\Big(\frac{e^{t\Delta}|\phi|}{\lambda}\Big)dx
\le
 \int_{B}e^{t\Delta}f\Big(\frac{|\phi|}{\lambda}\Big)dx
 \le
 \int_{B}f\Big(\frac{|\phi|}{\lambda}\Big)dx.
\]
This yields the desired estimate.

\begin{lemma}\label{ioku} Assume $1\leq p\leq 2$. There exists a positive constant $C$
such that
$$
\|e^{t\Delta}\phi\|_{L^f_\gamma(B)}\leq  C \, t^{-\frac
1p}\Big(\log\big(t^{-1}+ 1\big)\Big)^{-1/2}\|\phi\|_{L^p(B)} \ \
$$
for all $\phi\in L^p(B)$, $t>0$.
\end{lemma}
This lemma in the whole space $\mathbb{R}^n$ was proved in \cite[Lemma~2.2]{I}.
The same method works in $B_{\rho}$ since
we only need the $L^p-L^q$ estimate of the heat kernel which still holds in $B_{\rho}$.

\subsection{Lorentz spaces and heat kernel}\label{L}
 We present some regularizing properties of the heat
kernel in Lorentz spaces.  We recall the definition of Lorentz
spaces $L^{p,q}(B)$ on a ball $B\subset \mathbb{R}^2$. Let $\phi$
be a measurable function on $B$, which is finite almost everywhere.  We define
the distribution function
$$ \mu(\lambda,\phi)=|\{x\in B: |\phi(x)|>\lambda \}|, \ \ \
\lambda \geq 0.
$$
The decreasing rearrangement of $\phi$ is the function $\phi^*$
defined on $[0,\infty)$ by
$$ \phi^*(t)=\inf \{\lambda>0: \mu(\lambda,\phi)\leq t\}, \ \ \ t\geq 0.
$$
The Lorentz space $L^{p,q}(B)$, with $1\leq p<\infty$   consists
of all $\Phi$ measurable on $B$ and finite a.e. for which the
quantity
$$
\begin{aligned}
&
\|\Phi\|^*_{L^{p,q}(B)}=\Big(\int_0^{\infty}(t^{1/p}\Phi^*(t))^q\frac{dt}{t}\Big)^{1/q}
\ \ \ \ \ \
&&
{\rm when}\ \ 1\leq q<\infty,
\\
&
\|\Phi\|^*_{L^{p,\infty}(B)}= \sup_{t>0}t^{1/p}\Phi^*(t) \ \  \ \ \ \ \ \
&&
{\rm when}\ \
q=\infty
\end{aligned}
$$
is finite. In general, $ \|\cdot \|^*_{L^{p,q}(B)}$ is a
quasi-norm, but when $p>1$ it is possible to replace the
quasi-norm with a norm, which makes $L^{p,q}(B)$ a Banach space.
In the following we will denote by $ \|\cdot \|_{L^{p,q}(B)}$ this
norm (see \cite[Section~7.25]{AF}).

The Lorentz spaces can also be defined using Schwarz
symmetrization $\Phi^\sharp$ of $\Phi$, given by $\Phi^\sharp(x)
:= \Phi^*(\pi|x|^2)$; therefore  $\Phi\in L^{p,q}(B)$, $1\leq p<
\infty$, if and only if
\[
\begin{aligned}
&
\Big(\int_{B}\Big(|x|^{\frac{2}{p}}\Phi^{\sharp}(x)\Big)^q
        \frac{dx}{|x|^2}
    \Big)^{\frac1q}<\infty  \ \ \ \ \
    &&
    {\rm when}\ \ 1\leq q<\infty,
    \\
&
 \sup_{x\in B}|x|^{2/p}\Phi^{\sharp}(x)<\infty,\ \  \ \
\ \ \ \
&&
{\rm when}\ \ q=\infty.
\end{aligned}
\]

 \begin{lemma}\label{lemma:1.1} Let  $1\leq q<\infty$ and $1<p\leq r<\infty$.
 There exists a positive constant $C>0$ such that
  $$t^{1/p-1/r}\|{\rm e}^{t\Delta}\phi\|_{L^{r,q}(B)}\leq
C\|\phi\|_{L^{p,q}(B)}\ \ \ \ \ \ \ { for\ all\ }t>0.$$ Moreover
for $1<p<r<\infty$ and  for all $\phi\in L^{p,q}(B)$ we have
\begin{equation}\label{eq:3.1}
 \lim_{t\to
0}t^{1/p-1/r}\|{\rm e}^{t\Delta}\phi\|_{L^{r,q}(B)}=0.
\end{equation}

 \end{lemma}

\noindent {\bf Proof.} The first assertion in the lemma is proved
by the $L^p$-$L^q$ estimate of the heat kernel
(see~\cite[Proposition~48.4]{QS}) and real interpolation methods
(see~\cite[Theorem~5.3.2]{BL2}). The second assertion is a
consequence of the density of $C_0^{\infty}$ in $L^{p,q}(B)$ with
$1\le q <\infty$.

\subsection{Weak and classical solutions}\label{3.4}
We now present the notions of weak and classical solution for the
Cauchy problem \eqref{eq:2.1} with initial data $u_0\in
\exp L^2(B_\rho)$ where  $B_{\rho}$ is the ball centered at the
origin and of radius $\rho>0$. For the sake of simplicity  we will
omit the underlying space $B_\rho$.

\begin{definition}[Weak solution]\label{definition:4.1}$ $
\par \noindent
Let $u_0\in \exp L^2$ and $u\in L^\infty (0,T;\exp L^2)$ for some
$T\in(0,+\infty]$. We call $u$ a weak solution of the Cauchy
problem \eqref{eq:2.1} if $u$ satisfies the differential equation
$\partial_tu-\Delta u
    =f(u)$ in ${\cal D}'((0,T)\times B_{\rho})$ and
    $u(t)\to u_0$ in weak$^*$ topology as $t\to 0$.
\end{definition}
We recall that $u(t)\to u_0$ in weak$^*$ topology as $t\to 0$ if
and only if
\[
\lim_{t\to 0}\int_{B_\rho}\big(u(t,x) -u_0(x)\big)\, \psi(x) dx=0
\]
for every $\psi$ belonging to the predual space of $\exp L^2$. The
predual space of $\exp L^2$ is the Orlicz space defined by the
complementary function of $A(t)={e}^{t^2}-1$, denoted by
$\widetilde A(t)$. This complementary function is a convex
function such that $\widetilde A(t)\thicksim t^2$ as $t\to 0$ and
$\widetilde A(t)\thicksim t \log^{1/2} t$ as $t\to +\infty$.

 \begin{definition}[Classical  solution]\label{definition:4.2}$ $
\par
\noindent
 Let
$u_0\in \exp L^2$ and $u\in C((0,T],\exp L^2)\cap L^\infty_{loc} ((0,T),
L^\infty)$ for some $T\in(0,+\infty]$. We say that  the function
$u$ is an $\exp L^2$-classical solution of the Cauchy problem
\eqref{eq:2.1} in $(0,T]$ if
    $\left\|u(t)-{\rm e}^{t\Delta} u_0\right\|_{\exp L^2}\to 0$  as $t\to 0$,   $u$ is $C^1$ in
$t\in (0,T)$, $C^2$ in $x\in B_{\rho}$, continuous on
$\overline{B}_{\rho}$ and  $u$ is a classical solution
\eqref{eq:2.1} on $(0,T)\times B_{\rho}$.
\end{definition}
 We remark that  any $\exp L^2$-classical solution of the Cauchy
 problem \eqref{eq:2.1}  is also a weak solution. Indeed we have  that $u\in L^\infty
 (0,\varepsilon; \exp L^2)$ for some  $\varepsilon>0$ and this is a consequence
 of the inequality
  $$
 \|u(t)\|_{\exp L^2}\leq\left\|u(t)-{\rm e}^{t\Delta}
 u_0\right\|_{\exp L^2}+\left\|{\rm e}^{t\Delta}
 u_0\right\|_{\exp L^2}
 $$
and
 $$
 \left\|u(t)-{\rm e}^{t\Delta} u_0\right\|_{\exp L^2}\to 0, \ \ t\to 0.
 $$
Finally $u(t)\to u_0$ in the weak* topology as $t\to 0$ since $
{\rm e}^{t\Delta} u_0\to u_0$ in the weak* topology as $t\to 0$
and $u(t)-{\rm e}^{t\Delta}
 u_0\to 0$ in $\exp L^2$.

\section{Construction of a singular stationary solution}\label{sec:sing}
In this section we prove the existence of a radial singular
solution for the Dirichlet boundary value problem \eqref{eq:1.2}
in $B_{\rho} \subset \mathbb R^2$, for a well chosen $\rho>0$, by
using the shooting method (see \cite{C} and \cite{MTW}); that is,
we give the
\par \medskip \noindent
{\bf Proof of Proposition \ref{prop1}}.
\par \smallskip \noindent
Defining
\[
U(r) = \sqrt{-2\log r},
\]
one easily checks  that
$U$ solves
$$-U'' - \frac 1 r \, U' = \frac 1{U^3} \ e^{U^2},\quad 0<r<1.
$$
%
%
The solution $U$ was found by de Figueiredo-Ruf in \cite[p.\,653]{deFR}.

Let $f(s)$ as in \eqref{eq:ioku1.2}.  We want to continue the solution $U$ to a solution of
\begin{equation}\label{e2}
\left\{ \begin{array}{ll}
\displaystyle - u'' - \frac 1 r  u' = f(u)  \ &\hbox{ in }\ (0,\rho), \vspace{0.2cm} \\
\quad u(\rho) = 0, \vspace{0.2cm}
\\
\quad u(r)>0\ &\hbox{ in }\ (0,\rho),
\end{array} \right.
\end{equation}
where $\rho$ will be determined later.

\noindent Note that the solution $U(r) = \sqrt{-2\log r}$
satisfies
$$
U(r) \ge \sqrt{\frac 52}
 \
 \iff \
r \le \frac 1{e^{5/4}}.
$$
Let us consider the following equation
\begin{equation}\label{e3}
\left\{ \begin{array}{ll} - v'' - \frac 1 r \ v' = \alpha\, v^2 \ , \ r\geq \frac 1{e^{5/4}}\, , \vspace{0.2cm} \\
v\Big(\frac 1{e^{5/4}}\Big) = \sqrt {\frac 52}, \vspace{0.2cm} \\
v'\Big(\frac 1{e^{5/4}}\Big) = U'\Big(\frac 1{e^{5/4}}\Big)=
-\frac {e^{5/4}}{\sqrt {5/2}}.
\end{array} \right.
\end{equation}
We now prove that there exists a first zero $\rho > \frac
1{e^{5/4}}$ of the solution $v(r)$ of the problem~(\ref{e3}) by using
a shooting method and a contradiction argument.

%
%
By contradiction, assume that $ v(r) > 0$,  for all $r > \frac
1{e^{5/4}}$. Then $v'(r) < 0$, for all $r > \frac 1{e^{5/4}}$; if
not, there would exist $r_0$ with $v'(r_0) = 0$ and $v''(r_0) \ge
0$, but then $-v''(r_0) = \alpha v^2(r_0) > 0$, which is
impossible.
It follows from the above argument that
$v(r)$ has a limit $L \ge 0$, as $r
\to \infty$. We first show that $L = 0$. Indeed, consider the energy
$$ E(v,r) := \frac 1 2 |v'(r)|^2 + \frac \alpha 3 v(r)^3.
$$
Multiplying the equation of \eqref{e3}
by $v'(r)$, we obtain
$$
-v''(r) v'(r) - \frac1 r |v'(r)|^2 = \alpha v(r)^2\, v'(r)
$$
and so it follows
$$\frac d{dr} E(v,r) = v'(r)\, v''(r) + \alpha v(r)^2 \, v'(r) = -\frac 1 r\, |v'(r)|^2.
$$
This yields that $E(v,r)$ is decreasing, and hence
$$
|v'(r)|^2 \le 2\, E\Big(v,\frac 1 {e^{5/4}}\Big).
$$
Then, using again the equation of \eqref{e3}, we conclude for $r \to \infty$
$$
-v''(r) - \frac 1 r\, v'(r) = \alpha \, v(r)^2 \to \alpha\, L^2
$$
that
$$v''(r) \to -\alpha\,L^2,
$$
from which
we obtain $L = 0$.
We now derive a contradiction by using $L=0$. Observe that
$$\left(rv'(r) - \frac 1{e^{5/4}}\, v'\left(\frac 1{e^{5/4}}\right)\right)' = v'(r) + rv''(r) = -r\alpha \, v(r)^2
$$
and hence
\begin{equation}\label{e4}
r\, v'(r) - \frac 1{e^{5/4}}\, v'\big(\frac 1{e^{5/4}}\big) =
-\int_{1/e^{5/4}}^r s\, \alpha \, v(s)^2 ds.
\end{equation}
Therefore
$$
\begin{aligned}
 -r\, v'(r)
&= \int_{1/e^{5/4}}^r s\,\alpha\, v(s)^2 ds + \sqrt{\frac 25}\vspace{0.2cm} \\
& \ge \alpha v(r)^2 \int_{1/e^{5/4}}^r sds + \sqrt{\frac 25} \vspace{0.2cm} \\
& > \alpha v(r)^2 \, \frac {r^2} 2.\
\end{aligned}
$$
This implies that
$\ds \frac 1{v(r)} - \frac{\alpha\, r^2}4$\  is increasing.
Thus
$$
\ds \frac 1{v(r)} - \frac{\alpha\, r^2}4 > \frac 1{v\left(\frac
1{e^{5/4}}\right)}- \frac{\alpha}{e^{5/2}\, 4} = \sqrt{\frac 25} -
\frac 14\left(\frac 25\right)^{5/2}    > 0,
$$
which yields
$$
\frac 4\alpha\, r^{-2} > v(r)
$$
and
$$
\int_{1/e^{5/4}}^\infty r\, v(r)^2 dr\le \int_{1/e^{5/4}}^\infty
r\, \frac {16}{\alpha^2}\, r^{-4}dr < \infty.
$$
 It follows from (\ref{e4}) that there exists $A>0$ such that
$$
r\, v'(r) = \frac 1{e^{5/4}}\, v'\left(\frac 1{e^{5/4}}\right) -
\alpha\int_{1/e^{5/4}}^r s\, v(s)^2ds \to -A < 0.
$$
Hence
$$ v(r) = \int_{1/e^{5/4}}^r v'(s)ds + v\left(\frac1{e^{5/4}}\right) \le C\int_{1/e^{5/4}}^r -\frac As\, ds
\le  -AC(\log s) {\big|}_{1/e^{5/4}}^r \to -\infty \ \ \ {\rm as\
}r\to +\infty.
$$
This yields a contradiction, and hence there must exist a first
zero $\rho$ for $v(r)$.

\noindent By the above argument, we see that
$$
w(r) := \left\{ \begin{array}{ll} U(r), \ 0 < r < \frac 1{e^{5/4}}, \vspace{0.2cm} \\
v(r), \ \frac 1{e^{5/4}} \le r \le \rho,
\end{array}
\right.
$$
satisfies the equation \eqref{e2}.
In the following we define  $$\widetilde u(x)=w(|x|)= \left\{ \begin{array}{ll} U(|x|), \ 0 < |x| < \frac 1{e^{5/4}},\vspace{0.2cm} \\
v(|x|), \ \frac 1{e^{5/4}} \le |x| \le \rho.
\end{array}
\right.$$
 We stress that $\widetilde u$ belongs to
$C^2(B_{\rho}\setminus\{0\})\cap
C(\overline{B_{\rho}}\setminus\{0\})$, $\widetilde u(x)=0$ on
$|x|=\rho$ and
$$\widetilde u(x)= \sqrt{-2\log|x|}, \ \ \ \ |x|\leq \frac 1{e^{5/4}}
$$
and it is a classical solution of the elliptic equation on
$B_{\rho}\setminus\{0\}$.
%
%
\par \smallskip
It remains to prove
that the solution $\widetilde{u}$ satisfies the elliptic equation
in the sense of distributions in $B_{\rho}$.
 We use similar arguments as in \cite{BL1}, page 265
and in \cite{NS}, pages 261-262. Let $\varphi$ be a
$ C^{\infty}$ function with compact support in
$B_{\rho}$. We prove that
$$
\int_{B_{\rho}} \widetilde{u}\ \Delta \varphi+\ f(\widetilde{u})\ \varphi\
dx=0.
$$
Indeed let $\Phi(r)$ be a $C^{\infty}(\mathbb{R})$
function, $0\leq \Phi(r)\leq 1$ such that
\begin{equation*}
\Phi (r)=\left\{\begin{split} 1\ \ \ &\ \text{if}\ \ r<1/2,\\
0\ \ \ &\ \text{if}\ \ r\geq 1,\\
\end{split}\right.
\end{equation*}
and
$\Phi_\varepsilon(|x|)=\Phi\left(\frac{\log|x|}{\log\varepsilon}\right)$
for any $x\neq 0$  (these cut-off functions are the
same as those used in \cite{BL1}).  By a direct computation  for
small $\varepsilon>0$ we get $\Phi_\varepsilon (|x|)=1$ for
$|x|>\sqrt \varepsilon$ and $\Phi_\varepsilon (|x|)=0$ for
$|x|\leq \varepsilon$ and  for  $x\neq 0$, we get
$\Phi_\varepsilon(|x|)\to 1$ for  $\varepsilon \to 0^+$. By the
Dominated Convergence Theorem, since $\widetilde u$ and
\begin{equation}
\label{eq:4.4c}
f(\widetilde u)=
\left\{
\begin{aligned}&
\frac{1}{|x|^2(-2\log|x|)^{3/2}}\ \
&&
\text{if}\ \ 0<|x|<
\frac 1{{\rm e^{5/4}}},
\\
&\alpha \, v^2(|x|)\ \ \ \ \ \  \ \ \ \ \ \ \
&&
\text{if}\ \ \frac 1{{\rm e^{5/4}}}\leq|x|<\rho\\
\end{aligned}
\right.
\end{equation}
belong to $L^1(B_{\rho})$, we have
\begin{equation*}
\begin{aligned}
&\int_{B_{\rho}} \widetilde u \, \Delta  \varphi+ f(\widetilde u)\, \varphi\
dx\\&=\lim_{\varepsilon \to 0^+}\int_{B_{\rho}}\Phi_\varepsilon
 \widetilde u  \, \Delta \varphi+\Phi_\varepsilon f(\widetilde u)\ \varphi\ dx
\\&=\lim_{\varepsilon \to 0^+}\int_{B_{\rho}}\Phi_\varepsilon\,
 \Delta \widetilde u\, \varphi\
dx+2\int_{B_{\rho}}\nabla\Phi_\varepsilon\cdot \nabla \widetilde u \
\varphi\ dx+\int_{B_{\rho}}\Delta\Phi_\varepsilon \,\widetilde u\,
\varphi\ dx
+\int_{B_{\rho}}\Phi_\varepsilon\, f(\widetilde u)\, \varphi\ dx.\\
\end{aligned}
\end{equation*}
Since $\widetilde u$ is a classical solution of the elliptic
equation in $B_\rho\setminus\{0\}$ we obtain
\begin{equation*}
\begin{aligned}
&\lim_{\varepsilon \to 0^+}\int_{B_\rho}\Phi_\varepsilon\, \widetilde u\,
\Delta \varphi+\Phi_\varepsilon f(\widetilde u)\ \varphi\ dx
\\&=\lim_{\varepsilon \to 0^+}2\int_{B_{\rho}}\nabla\Phi_\varepsilon\cdot
\nabla \widetilde u\ \varphi \
dx+\int_{B_{\rho}}\Delta\Phi_\varepsilon\, \widetilde u\, \varphi\
dx.\\
\end{aligned}
\end{equation*}
Since $$ \Delta \Phi_\varepsilon=\Phi^{''}\Big(\frac{\log r}{\log
\varepsilon}\Big)\frac{1}{r^2(\log\varepsilon)^2}$$
 we have
 $$
 \Big|\int_{B_{\rho}} \widetilde u \, \Delta \Phi_\varepsilon \, \varphi\ dx\Big|\leq
 \frac{C}{(\log\varepsilon)^2}\int_{\varepsilon}^{\sqrt{\varepsilon}}\frac{\sqrt{-2\log (r )}}{r}\ dr
 $$
 and
 $$
 \begin{aligned}
\displaystyle  \lim_{\varepsilon\to
0^+}\frac{\int_\varepsilon^{\sqrt{\varepsilon}}\frac{\sqrt{-2\log
r}}{r}\
 dr}{(\log\varepsilon)^2}
 =
 \lim_{\varepsilon \to 0}
 \frac{2\sqrt{2}-1}{
 3\sqrt{-\log \varepsilon}}=0.
\end{aligned}
$$
In a similar way
$$
 \Big|\int_{B_{\rho}} \nabla \widetilde u\cdot \nabla \Phi_\varepsilon\ \varphi \ dx \Big|\leq
 \frac{C}{(-\log\varepsilon)}\int_{\varepsilon}^{\sqrt{\varepsilon}}\frac1{r\sqrt{-2\log r }}\ dr
 $$
and
$$
 \lim_{\varepsilon \to 0^+} \frac{\int_{\varepsilon}^{\sqrt{\varepsilon}}\frac1{r\sqrt{-2\log r }}\
 dr}{(-\log\varepsilon)}
 =
\lim_{\varepsilon \to 0^+} \frac{\sqrt{2}-1}{\sqrt{-\log
\varepsilon}}=0.
$$
This proves that the function $\widetilde u$ satisfies the
equation~\eqref{eq:1.2} in the sense of distributions.


\section{Well-posedness result}
In this section we consider   the Cauchy problem \eqref{eq:2.1}
where the  initial datum $u_0(x)$ is a measurable function
satisfying
\begin{equation}\label{eq:ioku5.1}
\mu:= \sup_{x\in B_{\rho}} \frac{u_0^{\sharp}(x)}{\widetilde
u(x)}<1.
\end{equation}
A typical example of such  initial data is $u_0=\mu \widetilde u(x)$
for $0<\mu <1$.
\par\noindent

Recall that $\int_{B_{\rho}}f(\widetilde u)\ dx<+\infty$
by \eqref{eq:4.4c}, hence
one can choose $\gamma=\int_{B_{\rho}}f(\widetilde u)dx$. With
this choice of $\gamma$, we now prove the  well-posedness result
1) in Theorem \ref{th2.1}. Let  $\max\{\mu,\frac{1}{\sqrt{2}}\}<\mu_1<1$ and
consider the complete metric space
$M_{T,\,\mu_1}$ introduced in \eqref{eq:2.4aaa}.
 We prove that there exist a positive time $T=T(\mu_1) $ and
a unique function $u\in M_{T,\,\mu_1}$ which is a
weak solution of \eqref{eq:2.1}. \par\noindent First, we make the
following:
\begin{remark}\label{re:5.1}
The initial data satisfying \eqref{eq:ioku5.1} belong to $M_{T,\,
\mu_1}$. Indeed, the  definition of $\gamma$ and a
standard property of the rearrangement  yield that
\[
\|\widetilde u\|_{L^f_\gamma}=\inf\Big\{\lambda >0:
\int_{B_\rho}f\Big(\frac{\widetilde u}\lambda \Big)dx\leq \gamma
\Big\}=1 \ \ \text{and} \ \
\|u_0\|_{L^f_\gamma}=\|u_0^\sharp\|_{L^f_\gamma}\le \mu
\|\widetilde u\|_{L^f_\gamma}=\mu< \mu_1.
\]
\end{remark}
\par \bigskip \noindent
In order to prove Theorem \ref{th2.1}.\,1) we first remark that in
the space $M_{T,\,\mu_1}$ the differential equation \eqref{eq:2.1}
admits an equivalent integral formulation as stated in the
following proposition.
\begin{proposition}\label{equiv}
 Let $u_0$ be a measurable function such that $\mu=\sup_{x\in B_{\rho}}
\frac{u_0^{\sharp}(x)}{\widetilde u(x)}<\mu_1<1$,
$T\in (0,+\infty]$ and $u\in M_{T,\,\mu_1}$. The
following statements are equivalent:

\noindent i) $u$ is a weak solution of the equation \eqref{eq:2.1}
in $(0,T)\times B_\rho$;

\noindent ii) $u$ satisfies the integral equation
\begin{equation}\label{inteq}
u(t)=e^{t\Delta} u_0+\int _0^t e^{(t-s)\Delta}f(u(s)) ds \ \ \ \
{\rm on}\ (0,T)\times B_\rho
\end{equation}
 in the sense of
distributions and $u(t)\to u_0$ as $t\to 0$ in the weak$^*$
topology.
\end{proposition}
The key tool of the proof  of  Proposition \ref{equiv} is the
 following
lemma:
\begin{lemma}\label{lem}
Let $0<\mu_1<1$, $T\in (0,+\infty]$ and $u\in
M_{T,\mu_1}$. Then
$$
\sup_{t\in (0,T)}
\left\|f(u(t))\right\|_{L^{\frac1{{\mu_1}^2}}}\leq
\left(C(\beta,\alpha,\mu_1)
  \gamma\right)^{{\mu_1}^2}.
$$
\end{lemma}

\noindent {\bf
Proof of Lemma \ref{lem}.}
 Since
 $\|u(t)\|_{L^f_\gamma}\leq\mu_1$, for any $t\in
(0,T)$, we  control uniformly with respect to time the
${L^f_\gamma}$-norm of the nonlinearity:
\begin{equation*}
\begin{aligned}
 \left\|f(u(t))\right\|_{L^{\frac1{{\mu_1}^2}}}^{\frac1{{\mu_1}^2}}&=\int_{B_\rho}f(u(t))^{\frac1{{\mu_1}^2}}dx\\
 &=\int_{|u|\geq \beta}\frac{{\rm e}^{\left(\frac{u}{\mu_1}\right)^2}}{|u|^{\frac
 3{{\mu_1}^2}}}dx+\int_{|u|<
 \beta}\alpha^{\frac1{{\mu_1}^2}}|u|^{\frac
 2{{\mu_1}^2}}dx\\
 &\leq\int_{|u|\geq \beta} \beta^{3-\frac 3{{\mu_1}^2}}\ \frac{{\rm e}^{\left(\frac{u}{\mu_1}\right)^2}}{|u|^3}dx+
 \int_{|u|<
 \beta}\alpha^{\frac1{{\mu_1}^2}}
 \beta^{\frac2{{\mu_1}^2}-2}\ |u|^2dx\\
& \leq
C(\beta,\alpha,\mu_1)\int_{B_{\rho}}f\Big(\frac{u}{\mu_1}\Big)dx
\leq  C(\beta,\alpha,\mu_1)
  \gamma\ \ \ \ \ \ \ \ \ \
\end{aligned}
\end{equation*}
for all $t\in (0,T)$. This ends the proof of Lemma~\ref{lem}.
\bigskip

The proof of Proposition \ref{equiv}  relies on the previous lemma
and  follows the same lines as the proof of Proposition 2.1 in
\cite{FKRT}.

\noindent We are now in position to prove the first part of
Theorem \ref{th2.1}.

\medskip
\noindent  {\bf
Proof of Theorem \ref{th2.1}.\,1)}
\par \medskip \noindent Let us introduce the integral operator
\[
\Phi(u)(t)=e^{t\Delta} u_0+\int _0^t e^{(t-s)\Delta}f(u(s)) ds
\]
and  look for  a fixed point of $\Phi$ in
$M_{T,\,\mu_1}$.

\noindent  First we prove that  $\Phi$ maps the space
$M_{T,\,\mu_1}$ into itself for small $T$.  By
applying Lemma \ref{heatcont} to the linear term and Lemma
\ref{ioku} with $p=\frac 1{{\mu_1}^2}$ ($\frac
12<{\mu_1}^2<1$) we obtain
\begin{equation*}
 \begin{aligned}
 \|\Phi(u)(t)\|_{L^f_\gamma}&\leq\left\|e^{t\Delta} u_0\right\|_{L^f_\gamma}+\int _0^t\left\| e^{(t-s)\Delta}f(u(s))\right\|_{L^f_\gamma}
 ds\\&\leq\left\| u_0\right\|_{L^f_\gamma}+\int _0^t(t-s)^{-{\mu_1}^2}\left(\log\left((t-s)^{-1}+1\right)\right)^{-\frac 12}
 \left\|f(u(s))\right\|_{L^{\frac1{{{\mu_1}}^2}}}ds.
\end{aligned}
\end{equation*}
Since  $\|u_0\|_{L^f_\gamma}\leq\mu$ (Remark \ref{re:5.1}) and the
${L^f_\gamma}$-norm of the nonlinearity is controlled uniformly
with respect to time   (Lemma \ref{lem}) we get
\begin{equation*}
\|\Phi(u)(t)\|_{L^f_\gamma}\leq \mu+\left(C(\alpha,
\beta, {\mu_1})\gamma\right)^{{\mu_1}^2}
\int_0^t\left(t-s\right)^{-{{\mu_1}}^2}\left(\log\left((t-s)^{-1}+1\right)\right)^{-\frac12}\
ds.
\end{equation*}
Since ${\mu_1}^2<1$ and
\begin{equation*}
\int_0^t\left(t-s\right)^{-{{\mu_1}}^2}\left(\log\left((t-s)^{-1}+1\right)\right)^{-\frac12}\
ds\rightarrow 0\ \ \ \ \ {\rm for}\ \ \ t\rightarrow 0,
\end{equation*}
if $T$ is small enough we get for any $0<t<T$ that
$${\left(C(\beta,\alpha,{\mu_1})\gamma\right)^{{\mu_1}^2}}
\int_0^t\left(t-s\right)^{-{{\mu_1}}^2}\left(\log\left((t-s)^{-1}+1\right)\right)^{-\frac12}\
ds\leq {\mu_1}-\mu$$ and this proves that $\Phi(u)$
belongs to $M_{T,{\mu_1}}$.
\par
Let us now prove that the integral operator $\Phi$ is a
contraction from $M_{T,{\mu_1}}$ into itself. Let  $q$
be such that $1<q<\frac1{{{\mu_1}}^2}$. We have
\begin{equation*}
\begin{aligned}
 \|\Phi(u)(t)-\Phi(v)(t)\|_{L^f_\gamma}&\leq\int _0^t\left\|
 e^{(t-s)\Delta}\left(f(u(s))-f(v(s))\right)\right\|_{L^f_\gamma}
 ds\\
&\leq \int_0^t(t-s)^{-\frac
1q}\left(\log\left((t-s)^{-1}+1\right)\right)^{-\frac 12}
\big\|f(u(s))-f(v(s))\big\|_{L^q} ds.
\end{aligned}
\end{equation*}
Since
$$
\left|f(u)-f(v)\right|\leq |u-v|\left(|f'(u)|+|f'(v)|\right)
$$ we have
$$
 \|f(u)-f(v)\|_{L^q}\leq  \|u-v\|_{L^{\widetilde r}}\left(\|f'(u)\|_{L^r}+\|f'(v)\|_{L^{r}}\right)
 $$
 where
 $\frac1q=\frac 1 {\widetilde r}+\frac 1r$, for $\widetilde r$ large enough such that
 $q<r<\frac 1{{\mu_1}^2}$.
Since $B_{\rho}$ is bounded,
 the Orlicz space is embedded into the Lebesgue space
 $L^{\widetilde r}$ (with $1<\widetilde r<\infty$). Therefore we have
$\|u-v\|_{L^{\widetilde r}}\leq\|u-v\|_{L^f_\gamma}.$
 Now, since $r<\frac{1}{{\mu_1}^2}$
\begin{equation}
\left|f'(u)\right|^{r}=\left\{ \begin{aligned} &\left|2|u|-\frac
3{|u|}\right|^r\Big( \frac{{\rm e}^{u^2}}{|u|^3}\Big)^r\leq
C(\beta,\mu_1, r)\, f\Big(\frac
u{\mu_1}\Big),
&&
|u|\geq\beta,\\
&(2\alpha |u|)^r,
&&
|u|<\beta.
 \end{aligned}\right.
\end{equation}
Therefore, thanks to the embedding of the Orlicz space in any
Lebesgue space $L^{\widetilde r}$, for $1<\widetilde r<\infty$,
and  since $ \sup_{s\in (0,T)}\|u(s)\|_ {L^f_\gamma}\leq
\mu_1$ we have
\begin{equation*}
\begin{aligned}
\|f'(u)\|_{L^r}&\leq C(\beta,\mu_1,
r)\left(\int_{|u|\geq
\beta}f\left(\frac{u}{\mu_1}\right)dx\right)^{\frac 1 r}+\Big(\int_{|u|\leq \beta}(2\alpha|u|)^r dx\Big)^{\frac 1r}\\
&\leq C(\beta,\mu_1, r)\left(\gamma\right)^{\frac 1r}+C(\alpha,r)\|u\|_{L^f_\gamma}\\
&\leq C(\alpha, \beta,\mu_1, \gamma, r).
\end{aligned}
\end{equation*}
Thus it holds
\begin{equation*}
\|f(u)-f(v)\|_{L^q}\leq C\|u-v\|_{L^f_\gamma},
\end{equation*}
for a constant $C=C(\alpha,\beta,\mu_1,\gamma,r)$.
Therefore, for all $0<t<T$,
\begin{equation*}
\|\Phi(u(t))-\Phi(v(t))\|_{L^f_\gamma}\leq
C
\sup_{0<t<T}
\|u(t)-v(t)\|_{L^f_\gamma}
\int
_0^t(t-s)^{-\frac
1q}\left(\log\left((t-s)^{-1}+1\right)\right)^{-\frac 12} ds
\end{equation*}
and
\begin{equation}\label{eq:elide6.1}
\int _0^t(t-s)^{-\frac
1q}\left(\log\left((t-s)^{-1}+1\right)\right)^{-\frac 12} ds\to 0,
\ \ \ {\rm as}\ t\to 0
\end{equation}
since $1<q<\frac 1{{\mu_1}^2}.$ This ends the proof
of the contraction argument.\par \noindent  We next prove the
convergence to the initial data $\left\|u(t)-{\rm e}^{t\Delta}
u_0\right\|_{\exp L^2}\to 0$  as $t\to 0$.
By the
equivalence of $L^f_{\gamma}$ and $\exp L^2$
(Proposition~\ref{proposition:2.1}),
we prove $\displaystyle \lim_{t\to 0}\|u(t)-e^{t\Delta}u_0\|_{\exp
L^2}=0$. Take $q$ so that $1<q<1/{\mu_1}^2$.
Lemma~\ref{ioku} gives us that
\[
 \|u(t)-e^{t\Delta}u_0\|_{\exp L^2}
 \le
 \int_0^t (t-s)^{-\frac{1}{q}}\big(\log ((t-s)^{-1}+1)\big)^{-\frac12}\|f(u(s))\|_{L^q}ds.
\]
By \eqref{relor}, for any $s\in (0,t)$ we have
\[
\|f(u(s))\|_{L^q} \le C\Big(\int_{B_{\rho}}\big(e^{q{u}^2}-1\big)\
dx\Big)^{\frac{1}{q}} \le
C'\Big(\int_{B_{\rho}}f\Big(\frac{u}{\mu_1}\Big)\
dx\Big)^{\frac{1}{q}} \le C' \gamma^{\frac1q}
\]
for some $C,C'>0$. Thanks to \eqref{eq:elide6.1} this gives  $
\|u(t)-e^{t\Delta}u_0\|_{\exp L^2} \to 0 $ as $t\to 0$.

\noindent Moreover $u$
 belongs to
$L_{loc}^{\infty}(0,T;L^\infty)$ (and so it is a
$\exp L^2-$classical solution of \eqref{eq:2.1} on $(0,T)\times
{B_\rho}$). Indeed assume  $t>0$. We know that ${ e}^{t\Delta}u_0$
belongs to $L^\infty$. Moreover, thanks to Lemma \ref{lem} we get
\begin{equation*}
\Big\|\int_0^t e^{(t-s)\Delta}f(u(s)) \
ds\Big\|_{L^{\infty}}\leq\int_0^t
(t-s)^{-{\mu_1}^2}\left\|f(u(s))\right\|_{L^{\frac1{\mu_1}^2}}ds\leq
C \int_0^t (t-s)^{-{\mu_1}^2}\ ds<+\infty
\end{equation*}
for fixed $t>0$. Finally by standard arguments one may check that
the solution $u$ belongs to $C((0,T],\exp L^2)$.

\par \bigskip

\section{Existence and Non-uniqueness result}
 In this section we prove the existence of an $\exp L^2-$classical
solution for the Cauchy problem \eqref{eq:2.1} for any nonnegative
$u_0$ such that
\[
\mu=\sup_{x\in B_{\rho}}\frac{{u_0}^{\sharp}(x)}{\widetilde
u(x)}\le 1.
\]
\par \medskip \noindent
This will imply the non-uniqueness result.
\par \medskip \noindent
{\bf Non-uniqueness:} Since $\widetilde u^\sharp(|x|)=\widetilde
u(x)$, we obtain that for the initial datum $u_0=\widetilde u$
 and for any $\mu_2>1$ there exist a positive time
$T=T(u_0,\mu_2) $ and an $\exp L^2$-classical solution $u$  of the
system \eqref{eq:2.1} that belongs to $M_{T,\, \mu_2}$.
  We
recall that $\widetilde u$ is a stationary {\it singular} solution
of the system \eqref{eq:2.1}, it is not bounded and it belongs to
the class $M_{T,1}$.
Therefore the  Cauchy problem \eqref{eq:2.1} possesses for
$u_0=\widetilde u$  at least two weak solutions in $M_{T,\mu_2}$,
even though a weak solution is unique in $M_{T,\mu_1}$ for $\mu<
\mu_1<1$ as in Theorem~\ref{th2.1} 1).

\begin{corollary}\label{theo:nonu}
Assume that  $u_0=\widetilde u$. For any $\mu_2>1$  there exist a positive time
$T=T(u_0, \mu_2)$ and  at least two weak solutions on $(0,T)\times
B_{\rho}$ of the Cauchy problem \eqref{eq:2.1} in the space $M_{T,\,\mu_2}.$
\end{corollary}

\par \bigskip \noindent
{\bf  Proof of Theorem \ref{th2.1}.\,2)}
\par \medskip \noindent
The key idea of the proof is to introduce a suitable auxiliary
Cauchy problem with a well--chosen polynomial nonlinearity whose
solutions can be transformed to {\it supersolutions} of the Cauchy
problem \eqref{eq:2.1}. Then, applying Perron's monotone method it
is possible  to prove the existence of a solution of
\eqref{eq:2.1}. To derive the auxiliary equation we apply the
generalized Cole-Hopf transformation introduced in \cite{FI}.
Define
$$
F(u):=\int_u^{+\infty} \frac{1}{f(s)}ds, \ \ \ \
\ u> 0,
$$
where $f$ is the nonlinearity defined in \eqref{eq:ioku1.2}. Now
let $  v_0=\max\left\{ {\left(F(u_0)\right)^{-1/2},
\left(F(\beta)\right)^{-1/2}}\right\}$, where $\beta$ is as in
\eqref{eq:ioku1.2}. Since $(F(t))^{-1/2}$ is a nondecreasing
function we obtain
$$  v_0^{\sharp}(|x|)=\left\{\begin{aligned}&(F(u_0^{\sharp}(|x|))^{-1/2}\ \ \ \ \ \
&&
\text{if}\ \ u_0^{\sharp}(|x|)> \beta,\\
&\left(F(\beta)\right)^{-1/2}\ \ \ \ \
&&
\text{if}\ \
u_0^{\sharp}(|x|)\leq \beta,
\end{aligned}\right.
$$
for any $x\in B_{\rho}$. It follows from
the definition of $f$ in \eqref{eq:ioku1.2}
that
\begin{equation}\label{eq:6.1c}
F(s)=\int_s^{\infty} \frac{\eta^3}{e^{\eta^2}}d\eta=\frac{s^2+1}{2e^{s^2}}\qquad \text{for\ large\ }s.
\end{equation}
Combining \eqref{eq:6.1c} to
the assumption on $u_0$,
we have
$$
v_0^{\sharp}(|x|)\ \le \left\{\begin{aligned}&\frac{\sqrt
2}{|x|(1-2\log|x|)^{1/2}},\ \ \ \ \ \
&&|x|< \frac1{{\rm
e}^{5/4}},\\
&\left(F(\beta)\right)^{-1/2},\ \ \ \ \
&&\frac1{{\rm e}^{5/4}}\leq |x|\leq {\color{blue}\rho}.\\
\end{aligned}\right.
$$
Consider the Cauchy problem
\begin{equation}\label{eq:v}
\left\{
  \begin{split}
   &\partial_tv-\Delta v
    =\frac{v^3}{2}
    && \text{in}\ \ B_{\rho}, \ t>0,
    \\
   &v(t,x)
    =F(\beta)^{-\frac 12}
    && \text{on}\ \partial B_{\rho}, \ t>0,\\
    &v(0,x)=v_0(x).\ \ \ \ \
  \end{split}
 \right.
 \end{equation}
If the initial  datum of \eqref{eq:v} belongs to
$L^2$, one can obtain a time-local classical solution by standard
contraction mapping arguments developed by Weissler~\cite{W1} and
Brezis-Cazenave~\cite{BC}. We should remark that the initial
 datum $v_0$ belongs to
any Lorentz space $L^{2,q}$ with $q>2$  since

\[
 \begin{aligned}
 v_0 \in L^{2,q} \iff
\int_{B_\rho}\left(|x|v_0^{\sharp}(x)\right)^q \frac{dx}{|x|^2}
<\infty
 \end{aligned}
\]
 and this last inequality is implied by the finiteness of the
 integral
\[
 \int_{|x|<e^{-5/4}}\frac{dx}{|x|^2 (1-2\log |x|)^{q/2}} <\infty.
\]
 We remark that $ v_0$ might not belong to $L^2$, as is the
case for  $u_0=\widetilde u$. Hence we consider the problem
\eqref{eq:v} in Lorentz space and obtain the following existence
result by modifying the arguments in \cite{W1,BC}.
\begin{proposition}\label{ex} Let $2<q\leq 5$. There exists a positive time $T=T(v_0)$ and a
unique solution $v$ of the Cauchy problem \eqref{eq:v} such that
$v\in C([0,T],L^{2,q})$, $t^{3/10}v(t)\in C([0,T],L^{5})$ and
$\lim_{t\to 0}t^{3/10}\|v(t)\|_{L^5}=0$.
Moreover $v\in
L^{\infty}_{loc}((0,T),L^\infty)$
and it is   a classical solution of \eqref{eq:v} on $(0,T)\times
B_{\rho}$.
\end{proposition}
We prove this proposition in the Appendix.
%
%
%

\par \smallskip
 We now build a super-solution of the Cauchy problem  \eqref{eq:2.1} by using the solution of \eqref{eq:v}.
Let us define
\[ \bar {u}=F^{-1}(v^{-2})
\]
where $F^{-1}$ is the inverse function of $F$ and $v$ is the
solution constructed in Proposition~\ref{ex}. Then $\bar u$
belongs to
$L^{\infty}_{loc}((0,T),L^\infty)$
because $v$ belongs to
$L^{\infty}_{loc}((0,T),L^\infty)$
 and $F^{-1}$ is a non-increasing function. Moreover,   $\bar{u}\ge 
F^{-1}(F(\beta))=\beta$, since $v(x,t)\geq(F(\beta))^{-1/2}$.
 Now by a direct computation we obtain
 $$
 \begin{aligned}
 \partial_t\bar{u}-\Delta \bar{u}-f(\bar{u})
 =4f(\bar{u})\, v^{-4}|\nabla v|^2 \left(\frac 32-f'(\bar{u})F(\bar{u})\right)\geq 0
\end{aligned}
$$
 since $f'(\bar{u})F(\bar{u})\leq 1$ for  any $\bar{u}\geq \beta$.
 Therefore,
 \begin{equation}\label{eq:dif}
 \partial_t\bar{u}\geq \Delta \bar{u}+f(\bar{u})
 \end{equation}
on $(0,T)\times B_{\rho}$. Moreover
$\bar{u}(0,x)=F^{-1}(v_0(x)^{-2})\geq u_0(x)$. Therefore, the
transformed function $\bar u$ is a supersolution of the original
problem \eqref{eq:2.1}. Applying Perron's monotone method, we
obtain a classical solution of the problem \eqref{eq:2.1} and of the corresponding integral equation
\eqref{inteq}\ (for more details, see
\cite[Proposition~2.1,  Lemma~2.3, Remark 6, (1)]{FI}).
%

We prove now the convergence of $u$
to the initial data, as $t \to 0$. We apply the following result.
\begin{lemma}[{\cite[Lemma~3.1]{FI}}]
Let $g(t)=f(F^{-1}(t))$. Assume that there exists  some $s_1>0$
such that
$$
f'(s)F(s)\leq 1 \ \ \ \ {\rm for\ all\ } s\geq s_1.
$$
Then there exists a constant $C$ such that
$
g(t)\leq Ct^{-1}
$
for all  $t<t_0=F(s_1)$.
\end{lemma}

\bigskip
Since
$$
u(x,t)\leq \bar{u}(x,t)=F^{-1}(v^{-2}(x,t)) \ \ \ \ {\rm and}\ \ \
\ v^{-2}(x,t)\leq F(\beta),
$$
by applying the previous
lemma
 we get
$$
\begin{aligned}
\left|u(t)-{\rm e}^{t\Delta}u_0\right|&=\int_0^t{\rm
e}^{(t-s)\Delta}f(u(s))ds\\&\leq\int_0^t{\rm
e}^{(t-s)\Delta}f(F^{-1}(v^{-2}(s))ds\\& \leq C\int_0^t{\rm
e}^{(t-s)\Delta}v^2(s)ds.
\end{aligned}
$$
Therefore
\begin{equation}\label{eq:6.3c}
\begin{aligned}
\left\|u(t)-{\rm e}^{t\Delta}u_0\right\|_{L^{\infty}}
&\leq C\, \Big\| \int_0^t{\rm e}^{(t-s)\Delta}v^2(s)ds\Big\|_{L^{\infty}}\vspace{0.2cm}\\
&\leq
C\int_0^t\frac{1}{(t-s)^{2/5}s^{3/5}}\, ds\
\Big(\sup_{0<s<t}s^{3/10}\|v(s)\|_{L^{5}}\Big)^2\\
& \leq
C\Big(\sup_{0<s<t}s^{3/10}\|v(s)\|_{L^{5}}\Big)^2
\end{aligned}
\end{equation}
and $\lim_{t\to 0}\sup_{0<s<t}s^{3/10}\|v(s)\|_{L^{5}}=0$. This
implies that $\|u(t)-{\rm e}^{t\Delta}u_0\|_{\exp L^2}\to 0$.

Furthermore, for any $\mu_2>\mu$, using also Lemma
\ref{heatcont}, we have
$$
\begin{aligned}
\sup_{0<t<T}\left\|u(t)\right\|_{L^f_\gamma}&\leq\sup_{0<t<T}\left\|{\rm
e}^{t\Delta}u_0\right\|_{L^f_\gamma}+\sup_{0<t<T}\left\|u(t)-{\rm
e}^{t\Delta}u_0\right\|_{L^f_\gamma}\\
&\leq\left\|u_0\right\|_{L^f_\gamma}+\sup_{0<t<T}\left\|u(t)-{\rm
e}^{t\Delta}u_0\right\|_{L^\infty}\\
&\leq\mu +\sup_{0<t<T}\left\|u(t)-{\rm
e}^{t\Delta}u_0\right\|_{L^\infty}\\
&\leq\mu_2
\end{aligned}
$$
for $T$ sufficiently small.
Hence $u\in M_{T,\mu_2}$.
Note that
\eqref{eq:6.3c} also implies that
$u\in L^{\infty}_{loc}(0,T; L^\infty)$. Finally, by standard
arguments one may check that  the solution $u$ belongs to
$C((0,T],\exp L^2)$.

\section {Non-existence result}
In this section we prove the  non-existence result
for $u_0 = \mu \, \widetilde u$ with $\mu > 1$, i.e. Theorem
\ref{th2.1}.3). We  start by stating the following:
\begin{proposition}\label{est}
Let $f$ be a $C^2$, positive, increasing, convex function in $(0,\infty)$ such that
$F(s):=\int_s^{\infty}\frac{1}{f(\eta)}d\eta <\infty$ for all $s>0$.
Let $u_0:B_{\rho}\to[0,\infty]$ and $u:B_{\rho}\times [0,T]\to
[0,\infty]$ be measurable functions satisfying
\begin{equation}\label{eqin}
u(t)\geq e^{t\Delta}u_0+\int_0^t e^{(t-s)\Delta}f(u(s)) ds \ \ \ \
\text{ \ a.e.\ in} \ \ \ B_{\rho}\times (0,T).
\end{equation}
Assume that $u(x,t)<\infty$ for a.e. $(x,t)\in B_{\rho}\times
(0,T)$.
Then there holds
\begin{equation}\label{upper}
\|e^{t\Delta}u_0\|_{L^{\infty}}\le F^{-1}(t)\ \ \ \ \ \ \text{
for\ all}\ \ t\in (0,T].
\end{equation}
\end{proposition}
{\bf Proof}.\ \
This proposition is essentially proved in \cite[Lemma~4.1]{FI} by applying the argument developed in Fujita~\cite[Theorem~2.2]{Fujita}
and Weissler \cite[Theorem~1]{W3}. Here we give a sketch of the proof for the reader's convenience.

Fix $\tau\in (0,T]$ and $t\in (0,\tau)$.
Applying $e^{(\tau-t)\Delta}$ to  \eqref{eqin}, we have by Fubini's theorem that
\[
 e^{(\tau -t)\Delta}u(t)\ge e^{\tau \Delta}u_0 + \int_0^t e^{(\tau -s)\Delta}f(u(s))ds
\]
for all $t\in (0,\tau)$.
Since $f$ is convex, one can apply Jensen's inequality to obtain
\begin{equation}\label{eq:2.1aaa}
 e^{(\tau-t)\Delta}u(t)\ge e^{\tau \Delta}u_0 + \int_0^t f\left(e^{(\tau-s)\Delta}u(s)\right)ds.
\end{equation}
Define
$\displaystyle
H(x,t):= e^{\tau \Delta}u_0 + \int_0^t f\left(e^{(\tau-s)\Delta}u(s)\right)ds.
$
Then we have
\[
-\frac{\partial}{\partial t}\left[F(H(x,t))\right]=\frac{\frac{\partial H}{\partial t}(x,t)}{f(H(x,t))}\ge 1.
\]
This yields
\[
 -F(H(x,t))+F(H(x,0))\ge t.
\]
Since $F(H(x,t))\ge 0$ and $H(x,0)=e^{\tau \Delta}u_0$,
there holds
\[
 e^{\tau \Delta}u_0 \le F^{-1}(t)
\]
for all $t\in (0,\tau)$.
Taking $t\uparrow \tau$ and the supremum on $x\in B_{\rho}$,
we obtain the desired estimate.

\begin{corollary}
\label{Corollary:7.1}
\noindent
Let $f$ be the function defined in \eqref{eq:ioku1.2}.
Assume that $u_0$ and $u$ satisfy the same conditions as in Proposition~\ref{est}.
Then there holds
\begin{equation}\label{ioku1015}
\|e^{t\Delta}u_0\|_{L^{\infty}} \le \left(-\log t\right)^{\frac{1}{2}} +1 \quad \text{for\ small\ }t>0.
\end{equation}
\end{corollary}
\noindent
{\bf Proof}.\ \
By
\eqref{eq:6.1c},
we have
\[
\begin{aligned}
\lim_{t\to 0}
\left[
{F^{-1}(t)}-{(-\log t)^{\frac12}}
\right]
&
=
\lim_{s\to \infty}
\left[
s-\left(\log \frac{1}{F(s)}\right)^{\frac12}
\right]
\\
&=
\lim_{s\to \infty}
\left[
s-\left(s^2+\log \frac{2}{s^2+1}\right)^{\frac12}
\right]
=0.
\end{aligned}
\]
Hence there holds
\[
F^{-1}(t)
\le
\left(-\log t\right)^{\frac{1}{2}} +1 \quad \text{for\ small\ }t>0.
\]
This and Proposition~\ref{est} yield the conclusion.

\par \medskip

Now we are in the  position to prove
Theorem~\ref{th2.1}.3).
\par \bigskip  \noindent
{\bf
Proof of Theorem~\ref{th2.1}.3)}
\par \medskip  \noindent
Assume that there exists a non-negative  $\exp L^2$-classical
solution of \eqref{eq:2.1} with $u_0=\mu \widetilde u,\ \mu>1$.
For any $t>0$, $s>0$, $t+s<T$ we have
\[
u(t+s)\geq {e}^{t\Delta}u(s).
\]
For $s\to 0$ we get
\begin{equation}\label{eq:7.1}
u(t)\geq {e}^{t\Delta}u_0
\end{equation}
thanks to the definition of $\exp L^2$-classical solution and  the
weak$^*$ convergence of $u(s)\to u_0$ as $s\to 0$. Since $u$ is an
$\exp L^2$-classical solution for any $0<\tau<t<T$ we have
\begin{equation}\label{eq:7.2}
u(t)=e^{(t-\tau)\Delta}u(\tau)+\int_\tau^te^{(t-s)\Delta}f(u(s))
ds.
\end{equation}
 Thanks to \eqref{eq:7.1} and \eqref{eq:7.2} we get
\begin{equation*}
u(t)\geq e^{t\Delta}u_0+\int_\tau^te^{(t-s)\Delta}f(u(s)) ds,
\end{equation*}
and for $\tau\to 0$ by monotone convergence theorem we have:
\begin{equation*}
u(t)\geq e^{t\Delta}u_0+\int_0^te^{(t-s)\Delta}f(u(s)) ds.
\end{equation*}
 Therefore applying
%
Corollary~\ref{Corollary:7.1},
we get that  $u$ satisfies
\eqref{ioku1015}.
We now
prove an estimate of $\left\|{e}^{t\Delta}u_0\right\|_{\infty}$
from below  which is in  contradiction with \eqref{ioku1015}. Remark
that
\begin{equation*}
\begin{aligned}
\left\|{e}^{t\Delta}u_0\right\|_{\infty}&\geq\int_{B_{\rho}(0)}G(0,y,t)\mu\,
\widetilde u(y)\ dy\\&\geq\int_{B_{r}(0)}G(0,y,t)\mu
\sqrt{-2\log|y|}\ dy,
\end{aligned}
 \end{equation*}
where $r={\frac{1}{e^{5/4}}}$.
 Let us denote by $d=\rho-r$. It is possible to bound on the ball $B_r(0)$ the Dirichlet heat kernel $G$ associated
  to the ball $B_{\rho}$ from below  by the heat
 kernel for $\mathbb{R}^2$ (see \cite{VB}):
\begin{equation*}
G(0,y,t)\geq H(d,t)\frac{e^{-|y|^2/4t}}{4\pi t},
\end{equation*}
where
$$
H(d,t)=1-{e^{-d^2/t}\left(2+4\frac{d^2}{t}\right)}.
$$
Therefore
\begin{equation*}
\begin{aligned}
\left\|{e}^{t\Delta}u_0\right\|_{\infty}&\geq\int_{B_{r}(0)}G(0,y,t)\mu
\sqrt{-2\log|y|}\ dy\\
&\geq \int_{B_{r}(0)}H(d,t)\frac{e^{-|y|^2/4t}}{4\pi t}\mu
\sqrt{-2\log|y|}\ dy\\
&\geq H(d,t)\int_{|z|\leq rt^{-1/2}}\frac{e^{-|z|^2/4}}{4\pi
}\mu
\sqrt{-\log t-2\log|z|}\ dz,\\
\end{aligned}
 \end{equation*}
 where in the last inequality we replace $y=\sqrt t\, z$.
 For $a<1/2$  and for small values of $t$ we obtain
\begin{equation*}
\begin{aligned}
\int_{|z|\leq rt^{-1/2}}\hspace{-0.15cm}\frac{e^{-|z|^2/4}}{4\pi
}\mu \sqrt{-\log t-2\log|z|}dz&\geq\int_{|z|\leq
rt^{-a}}\hspace{-0.15cm}\frac{e^{-|z|^2/4}}{4\pi }\mu \sqrt{-\log t-2\log|z|}dz \\
&\geq \mu \sqrt{-\log t+2a\log t-2\log r}\int_{|z|\leq rt^{-a}}\hspace{-0.15cm}\frac{e^{-\frac{|z|^2}{4}}}{4\pi }dz \\
&\geq \mu\, \sqrt{1-2a} \sqrt{-\log t} \ (1-\varepsilon)\\
\end{aligned}
\end{equation*}
for some $\varepsilon > 0$, since $\int_{|z|\leq
rt^{-a}}\frac{e^{-|z|^2/4}}{4\pi }dz\to 1 $ for $t\to 0^+$. Since
also $H(d,t) \to 1$ as $t \to 0^+$, we get
\begin{equation*}
\begin{aligned}
\left\|{e}^{t\Delta}u_0\right\|_{\infty} &\geq
\mu \ H(d,t)(1-\varepsilon)\sqrt{1-2a} \sqrt{-\log t}\\
&\geq \mu \ (1-\varepsilon)^2 \sqrt{1-2a} \ \sqrt{\log
\frac 1 t}
\end{aligned}
\end{equation*}
Thus, for fixed $\mu > 1$ we can choose $\varepsilon
> 0$ small and $a$ near $0$ such that
$$\mu (1-\varepsilon)^2 \sqrt{1-2a} \ge  1 + \delta
$$
for some $\delta > 0$.
 This contradicts \eqref{ioku1015} in the limit $t \to 0$.

\section{Appendix}
Proposition~\ref{ex} can be proved by a
 modification of the
standard contraction mapping argument developed by
Weissler~\cite{W1} and Brezis-Cazenave~\cite{BC} to the  framework of
Lorentz spaces. We include it for the reader's convenience.
\par \bigskip \noindent
\noindent {\bf Proof of Proposition~\ref{ex}.} We look for a
solution $v=\bar v+F(\beta)^{-\frac 12}$ where $\bar v$ is a
solution of the following Cauchy problem with Dirichlet boundary
condition:
 \begin{equation}\label{eq:vbar}
\left\{
  \begin{split}
   &\partial_t\bar v-\Delta \bar v
    =\frac{\left(\bar v+F(\beta)^{-\frac 12}\right)^3}{2}
    && \text{in}\ \ B_{\rho}(0), \ t>0,
    \\
   &\bar v(t,x)
    =0
    && \text{on}\ \partial B_{\rho}(0), \ t>0,\\
    &\bar v(0,x)=\bar v_0(x)
    && \text{in}\ \ B_{\rho}(0),\ \ \ \ \
  \end{split}
 \right.
 \end{equation}
 where $\bar v_0(x)=v_0(x)-F(\beta)^{-\frac 12}$.
 We prove that there exists a  solution $\bar v$ of the equation
 \eqref{eq:vbar} belonging to the space
$$
E_{\delta, M,T}=\left\{
w\in L^{\infty}(0,T;
L^{2,q}):
\begin{aligned}
&\sup_{t\in(0,T)}\|w(t)\|_{L^{2,q}} \leq M+1,\\
&\sup_{t\in(0,T)}t^{3/10}\|w(t)\|_{L^5} \leq
\delta\\\end{aligned}\right\}$$ where $M\geq \sup_{t\in
(0,\infty)} \|{\rm e}^{t\Delta}\bar v_0\|_{L^{2,q}}$ and $\delta $
and $T$ are well-chosen positive constants.

\noindent  Let us first remark that  the space $E_{\delta, M,T}$
endowed with the metric
$$d(v,w)=\sup_{t\in
(0,T)}t^{3/10}\|v(t)-w(t)\|_{L^5}$$ is a nonempty complete metric
space. Let us denote  $D=F(\beta)^{-\frac 12}$ and consider the
integral operator
$$
G(w)(t)={\rm e}^{t\Delta}{\bar v}_0+\frac 12 \int_0^t{\rm
e}^{(t-s)\Delta}\left(w(s)+D\right)^3 \ ds.$$
 We prove  that for
some well-chosen positive constants $T$ and $\delta$ the operator
$G$ maps the space $E_{\delta, M,T}$ into itself and it is a
contraction. Indeed let $w\in E_{\delta, M,T}$; by the smoothing
effect of the heat semigroup established in Lemma \ref{lemma:1.1},
$e^{t\Delta }D\leq D$ for any positive constant $D$, and thanks to
the inequality $\left|w+D\right|^3\leq 4\left(|w|^3+D^3\right)$,
for $t\in (0,T)$, we have
\begin{equation*}
\begin{aligned}
t^\frac{3}{10} \|G(w)(t)\|_{L^5} & \leq t^\frac{3}{10}\|{\rm e}^{t\Delta}{\bar
v}_0\|_{L^5}+2t^\frac{3}{10}\int_0^t\left\|{\rm
e}^{(t-s)\Delta}\left(|w(s)|^3+D^3\right) \right\|_{L^5}ds\\
&\leq t^\frac{3}{10}\|{\rm e}^{t\Delta}{\bar
v}_0\|_{L^5}
+\int_0^t\frac{Ct^\frac{3}{10}}{(t-s)^{2/5}s^{9/10}}ds
\left(\sup_{0<s<t}s^\frac{3}{10}\|w(s)\|_{L^5}\right)^3+C t^{\frac {13}{10}}\\
&\leq t^\frac{3}{10}\|{\rm e}^{t\Delta}{\bar
v}_0\|_{L^5}+C_1\delta^3+C_2 t^{\frac {13}{10}}.
\end{aligned}
\end{equation*}

\noindent Therefore
\begin{equation*}
\sup_{{t\in (0,T)}}t^{3/10}\|G(w)(t)\|_{L^5}\leq\sup_{{t\in
(0,T)}}t^{3/10}\|{\rm e}^{t\Delta}{\bar
v}_0\|_{L^5}+C_1\delta^3+C_2 T^{\frac {13}{10}}.
\end{equation*}

\noindent Moreover,
since
$L^{2}\subset L^{2,q}\ (q>2)$
we obtain
\begin{equation}\label{app}
\begin{aligned}
\|G(w)(t)\|_{L^{2,q}}&\leq \|{\rm e}^{t\Delta}\bar
v_0\|_{L^{2,q}}+2\int_0^t \Big\|{\rm
e}^{(t-s)\Delta}(|w(s)^3+D^3) \Big\|_{L^{2}}ds \\
&\leq
M+ C\int_0^t\frac{1}{(t-s)^{1/10}s^{9/10}}ds\Big(\sup_{0<t<T}t^{3/10}\|w(t)\|_{L^{5}}\Big)^3+Ct\\
&\leq M+C_3\delta^3+C_4T.\\
\end{aligned}
\end{equation}
\noindent Therefore
\begin{equation*}
\sup_{{t\in (0,T)}}\|G(w)(t)\|_{L^{2,q}}\leq M+C_4\delta^3+C_3T.
\end{equation*}
 In a similar way, since $|(w+D)^3-(v+D)^3|\leq
C|w-v|(w^2+v^2+D^2) $, for any $v,w\in E_{\delta, M,T}$, we have
\begin{equation}\label{r3}
\begin{aligned}
t^\frac{3}{10}\|G(v)(t)-G(w)(t)\|_{L^5} &\leq C
t^\frac{3}{10}\int_0^t\left\|{\rm
e}^{(t-s)\Delta}|v(s)-w(s)|(v^2(s)+w^2(s)+D^2) \right\|_{L^5}ds\\
&\leq \sup_{t\in
(0,T)}t^\frac{3}{10}\|v(t)-w(t)\|_{L^5}\left(C_5\delta^2+C_6T\right).\\
\end{aligned}
\end{equation}

\noindent Thus we obtain
\begin{equation*}
\sup_{t\in (0,T)}t^{3/10}\|G(v)(t)-G(w)(t)\|_{L^5}\leq \sup_{t\in
(0,T)}t^{3/10}\|v(t)-w(t)\|_{L^5}\left(C_5\delta^2+C_6T\right).
\end{equation*}

\noindent  Therefore by choosing $\delta$ such that
$$
C_1\delta^2\leq\frac 12,\ \ \ C_3\delta^3\leq\frac 12,\ \ \  C_5
\delta^2\leq\frac 14$$ and  $T$ small enough such that
$$
\sup_{t\in (0,T)}t^{3/10}\|{e}^{t\Delta}\bar
v_0\|_{L^5}+C_2T^{\frac {13}{10}}\leq \frac \delta2,\ \
  C_4T\leq \frac 12, \ \ C_6T\leq\frac 14$$
 we obtain that $G$ maps $E_{\delta,M, T}$ into itself and it
is a contraction.
We remark that $\sup_{t\in
(0,T)}t^{3/10}\|{e}^{t\Delta}\bar v_0\|_{L^5}\to 0$ as $T\to 0$
since $\bar v_0\in L^{2,q}$, with $2<q\leq 5$, thanks to Lemma
\ref{lemma:1.1}.
%
Therefore, the integral equation
\begin{equation}\label{eqintv3}
w(s)={\rm e}^{t\Delta}{\bar v}_0+\frac 12\int_0^t{\rm
e}^{(t-s)\Delta}\left(w(s)+D\right)^3 \ ds
\end{equation}
 admits a
unique solution $\bar v$ in $E_{\delta,T,M}$.

We prove now that the fixed point $\bar v$ belongs to
$$
E=E_{\delta,M,T}\cap \Big\{w\in C((0,T],L^5): \lim_{t\to 0}t^{3/10}\|w(t)\|_5=0\Big\}.
$$
To this end, it is enough to prove that $\Phi$ is a map from $E$
to $E$, since this implies that the previous contraction mapping
argument works in $E$. It follows from $v_0\in L^{2,q}$ and
Lemma~\ref{lemma:1.1} that
$e^{t\Delta}\bar{v_0}=e^{t\Delta}(v_0-F(\beta)^{-\frac12})\in E$.
Fix $w\in E$; since $E\cap C([0,T],L^{\infty})$ is dense in $E$
with respect to the metric $d$, there exists a sequence $v_n\in
E\cap C([0,T],L^{\infty})$ such that $G(v_n)\in E$ and
$d(v_n,w)\to 0$ as $n\to \infty$.  By $\eqref{r3}$, we have
$d(G(v_n),G(w))\to 0$ as $n\to \infty$. This together with the
fact that $E$ is a complete metric space with respect to $d$
yields $G(w)\in E$. This proves that the fixed point $\bar v$
belongs to $E$. Furthermore, \eqref{app} and $\bar v \in E$ yield
\[
 \lim_{t\to 0} \|\bar v(t)-e^{t\Delta}\bar{v}_0\|_{L^{2,q}}=0.
\]

Finally we prove that $\bar v$ is a classical solution.  Since
$\bar v_0$ is nonnegative,  the solution $\bar v$ is also
nonnegative. Moreover, it belongs to $ L^{\infty}_{loc}(0,T;
L^{\infty})$ and it is a classical solution on $(0,T)\times
B_{\rho}$. Indeed
\begin{equation}
\begin{aligned}
\|\bar v(t)\|_{L^{\infty}}&\leq \|{\rm e}^{t\Delta}\bar
v_0\|_{L^{\infty}}+\frac 12\Big\|\int_0^t{\rm
e}^{(t-s)\Delta}(\bar v+D)^3 \ ds \Big\|_{L^{\infty}} \\
&\leq \|{\rm
e}^{t\Delta}\bar v_0\|_{L^{\infty}}+ C\int_0^t\frac{1}{(t-s)^{3/5}}\, \|\bar v^3(s)\|_{L^{5/3}}ds+tCD^3\\
&\leq t^{-1}\|
v_0\|_{L^{1}}+C\int_0^t\frac{1}{(t-s)^{3/5}s^{9/10}}\, (s^{3/10}\|\bar v(s)\|_{L^{5}})^3ds+tCD^3\\
&\leq t^{-1}\|v_0\|_{L^{1}}+C t^{-1/2}\big(\sup_{s\in
(0,t)}s^{3/10}\|\bar v(s)\|_{L^{5}}\big)^3+tCD^3.
\end{aligned}
\end{equation}
Therefore, for any $\epsilon >0$, $v\in L^{\infty}(\epsilon,T;
L^\infty)$ and $\bar v$ is a classical solution on $(0,T)\times
B_\rho(0)$.

By denoting $v(x,t)=\bar v(x,t)+D$, $D=F(\beta)^{-1/2}$, we obtain
a solution of the differential equation \eqref{eq:v}. The solution
$v$ of \eqref{eq:v} belongs to $C([0,T],L^{2,q})\cap
C((0,T],L^{5})$ and $\lim_{t\to 0}t^{3/10}\|v(t)\|_{L^5}=0$ and it
is bounded on any interval $(\varepsilon, T)$, for
$\varepsilon>0$. Moreover $v(x,t)\geq F(\beta)^{-1/2}$ for any
$(x,t)\in B_{\rho}\times (0,T).$




\par 
{\it Adresses:}\par \smallskip
Norisuke Ioku*
\par
{\small Graduate School of Science and Engineering
\par
Ehime University
\par
Matsuyama, Ehime 790-8577, Japan}
\par
Email: {\tt ioku@ehime-u.ac.jp}
\par \medskip
Bernhard Ruf
\par
{\small Dipartimento di Matematica,
Universit\`a di Milano
\par
via C. Saldini 50, Milano 20133, Italy
}
\par
Email: {\tt bernhard.ruf@unimi.it}
\par \medskip
Elide Terraneo
\par
{\small Dipartimento di Matematica, Universit\`a di Milano
\par
via C. Saldini 50, Milano 20133, Italy}
\par
Email: {\tt elide.terraneo@unimi.it}
\par \medskip \noindent
{\small * This author was partially supported by JSPS  Grant-in-Aid for Young Scientists B \#15K17575.}

%
\end{document}